\documentclass[11pt,leqno]{amsart}
\usepackage{amsmath,amsfonts,amsthm,amscd,amssymb,graphicx}
\usepackage{epstopdf}
\numberwithin{equation}{section}
\usepackage{fullpage}

\newcommand{\imply}[1]{\Rightarrow}

\newcommand{\Ea}{\ensuremath{\mathcal{E}_A}}
\renewcommand\d{\partial}

\renewcommand\d{\partial}

\renewcommand\o{\omega}

\def\eps{\varepsilon }
\def\e{\varepsilon}

\renewcommand\d{\partial}

\renewcommand\o{\omega}
\newcommand\R{\mathbb R}
\newcommand\C{\mathbb C}

\def\eps{\varepsilon}
\def\e{\varepsilon}

\newcommand\E{\mathcal{E}}

\newcommand\br{\begin{remark}}
\newcommand\er{\end{remark}}
\newcommand\brs{\begin{remarks}}
\newcommand\ers{\end{remarks}}
\newcommand\bp{\begin{pmatrix}}
\newcommand\ep{\end{pmatrix}}
\newcommand{\be}{\begin{equation}}
\newcommand{\ee}{\end{equation}}
\newcommand\ba{\begin{equation}\begin{aligned}}
\newcommand\ea{\end{aligned}\end{equation}}

\newcommand{\bap}{\begin{app}}
\newcommand{\eap}{\end{app}}
\newcommand{\begs}{\begin{exams}}
\newcommand{\eegs}{\end{exams}}
\newcommand{\beg}{\begin{example}}
\newcommand{\eeg}{\end{exaplem}}
\newcommand{\bpr}{\begin{proposition}}
\newcommand{\epr}{\end{proposition}}
\newcommand{\bt}{\begin{theorem}}
\newcommand{\et}{\end{theorem}}
\newcommand{\bc}{\begin{corollary}}
\newcommand{\ec}{\end{corollary}}
\newcommand{\bl}{\begin{lemma}}
\newcommand{\el}{\end{lemma}}
\newcommand{\bd}{\begin{definition}}
\newcommand{\ed}{\end{definition}}

\newcommand{\CalT}{\mathcal{T}}

\newcommand{\Id}{{\rm Id }}

\newtheorem{theorem}{Theorem}[section]
\newtheorem{proposition}[theorem]{Proposition}
\newtheorem{corollary}[theorem]{Corollary}
\newtheorem{lemma}[theorem]{Lemma}

\theoremstyle{remark}
\newtheorem{remark}[theorem]{Remark}
\newtheorem{remarks}[theorem]{Remarks}
\theoremstyle{definition}
\newtheorem{definition}[theorem]{Definition}
\newtheorem{assumption}[theorem]{Assumption}

\newtheorem{example}[theorem]{Example}

\newcommand\cA{{\mathcal { A}}}

\newcommand\cE{{\mathcal  E}}

\newcommand{\beq}{\begin{equation}}
\newcommand{\eeq}{\end{equation}}

\renewcommand{\sp}{{\ensuremath{{\scriptscriptstyle +}}}}
\newcommand{\sm}{{\ensuremath{{\scriptscriptstyle -}}}}

\title{
	Recent results on stability of planar detonations
}

\author{Kevin Zumbrun}
\address{Indiana University, Bloomington, IN 47405}
\email{kzumbrun@indiana.edu} 
\thanks{Research of K.Z. was partially supported
under NSF grants no. DMS-0300487 and DMS-0801745.}
\begin{document}

\begin{abstract}
We describe recent analytical and numerical results on stability and behavior of viscous and inviscid
detonation waves obtained by dynamical systems/Evans function techniques like those used to study 
shock and reaction diffusion waves.
In the first part, we give a broad description of viscous and inviscid results for 1D perturbations;
in the second, we focus on inviscid high-frequency stability in multi-D and 
associated questions in turning point theory/WKB expansion.
\end{abstract}

\date{\today}
\maketitle

\textbf{Dedicated to Guy M\'etivier on the occasion of his 65th birthday.}

\tableofcontents

In these notes, we describe some recent work on stability and behavior of detonation waves,
carried out from a point of view evolving from the study of viscous and inviscid shock and boundary layers
in, e.g., \cite{GZ,ZH,Br,ZS,Z1,MZ,GMWZ1,GMWZ2,HuZ,HLZ,HLyZ1,HLyZ2,PZ}.
This material was originally presented as a pair of 90-minute lectures 
at the INDAM conference {\it Nonlinear Optics and Fluid Mechanics}, given in Rome, September 14-18, 2015
in honor of the 65th birthday of Guy M\'etivier, and our treatment follows closely to the spirit and
format of the lectures.

The topic was chosen for interest of the honoree as almost the unique one studied by
the author on which he has not explicitly collaborated with M\'etivier; 
nonetheless, many of the ideas may be seen to be related to ideas and tools developed by and
with Guy in other contexts.
The material presented here was developed in joint work with
Blake Barker, Jeff Humperys, Olivier Lafitte, Greg Lyng, Reza Raoofi, Ben Texier, and Mark Williams.
We mention also the foundational work of Kris Jenssen together with Lyng and Williams \cite{JLW},
of which we make frequent use.

\section{Stability of viscous and inviscid detonation waves}\label{s:I}

In this first part, we survey a collection of theoretical and numerical results on {\it 1D stability of detonations}
obtained over the past 10-15 years via Evans function-based techniques like those used to study shock and reaction diffusion waves.  These include stability in the small heat-release and high-overdrive limits, 
rigorous characterization of 1D instability as ``galloping'' type Hopf bifurcation,  description of the inviscid (ZND) limit, 
and numerical computation of viscous (rNS) spectra revealing a new phenomenon of ``viscous hyperstabilization.''

\medskip

Two underlying questions we have in mind in this section are:
\medskip

$\bullet$ {\it What is the (physical or mathematical) role of viscosity in the theory?}

\medskip

$\bullet$
{\it What is our role in the theory?  That is, what can we usefully contribute by our new techniques?}

\subsection{Viscous and inviscid detonation waves}\label{I.I}

Consider a general abstract combustion model, expressed in 1D Lagrangian coordinates
\cite{Z1,LyZ1,LyZ2,LRTZ,TZ4}:
\ba \label{rNS}
u_t + f(u)_x &= \eps (B(u)u_{x})_x +kq\phi(u)z,\\
z_t &= \eps (C(u,z)z_x)_x -k\phi(u)z,
\ea
$u$, $f$, $q\in \mathbb{R}^n$, $B\in \mathbb{R}^{n\times n}$,
$z\in \R^r$, $k$, $C$, $\phi \in R^1$, and $k,\, \eps>0$.
Here, $u$ comprises gas-dynamical variables,
$z=$ mass fraction(s) of reactant(s), 
$\phi(u)=$ ``ignition function'', 
$q=$ heat release, $k=$ reaction rate,
and $\eps$ (typically small) scales coefficients of viscosity/heat conduction/species diffusion.

A {\it right-going detonation solution} consists of a traveling wave
\begin{equation}
(u,z)(x,t)=(\bar u, \bar z)(x-st),
\qquad \lim_{x\to \pm \infty} (u,z)(x)=(u_\pm,z_\pm),
\notag \end{equation}
$s>0$, with $z_-=0$ and $z_+=1$, moving to the right into the totally unburned region toward $x\to +\infty$
and leaving behind the totally burned region toward $x\to -\infty$.

\begin{example}\label{rnseg}
A standard example is the reactive Navier--Stokes/Euler system
\begin{equation}\label{rNSeg}
\left\{ \begin{aligned}
 \d_t \tau - \d_x u & = 0,\\
 \d_t u + \d_x p  & = \d_x (\nu\tau^{-1} \d_x u),\\
 \d_t E + \d_x(pu) & = 
\d_x\big(
\kappa \tau^{-1} \d_x T + \nu\tau^{-1} u \d_x u\big)
 +qk \phi(T) z ,\\
 \d_t z 
& = \d_x (d \tau^{-2} \d_x z) - k \phi(T) z ,\\
\end{aligned}\right.
\end{equation}
where $\tau>0$ denotes specific volume, $u$ velocity, 
$ E = e+ \frac{1}{2} u^2$ 
specific gas-dynamical energy, 
$e>0$ specific internal energy,  
and $0 \leq z \leq 1$ mass fraction of the reactant,
with ideal gas equation of state, single-species reaction, and Arrhenius-type ignition function,
\be\label{eoseg}
 p=\frac{\Gamma  e}{\tau}, \quad T=c^{-1} e, 
\quad \phi(T)=e^{\frac{-\mathcal{E}}{T}},
\ee
For $\nu,\kappa, d>0$, this represents the ``viscous'' (mixed hyperbolic--parabolic) {\it reactive Navier--Stokes} (rNS)
equations \cite{Ba,CF}, for 
$\nu,\kappa, d=0$, the ``inviscid'' (hyperbolic) reactive Euler, or {\it Zel'dovich--von Neumann--D\"oring} (ZND) 
equations \cite{Ze,vN1,vN2,D}.
These represent successive refinements of the earlier {\it Chapman--Jouget} (CJ) theory \cite{C,J1,J2}, 
in which both transport (diffusion) and reaction processes are taken to occur instantaneously, across an ideal
shock-like discontinuity.
\end{example}

\subsubsection{Inviscid (ZND) Profiles}\label{s:param}
(Following \cite{Z2})
In case $\nu,\kappa,d=0$, $r=1$, we may explicitly solve the profile equation for \eqref{rnseg}--\eqref{eoseg}.
By the invariances of \eqref{rNSeg}--\eqref{eoseg}, we may take without loss of generality $\tau_+=1$, $u_+=0$, $s=1$, 
and $0 \leq e_+ \leq \frac{1}{\Gamma(\Gamma+1)}$, with  $\Gamma>0$, $\mathcal{E}>0$, $0\leq q \leq q_{cj}=
\frac{ (\Gamma+1)^2(\Gamma e_+ + 1)^2- \Gamma (\Gamma+2) ( 1+2(\Gamma +1)e_+ ) }
{2\Gamma(\Gamma+2)}$,
yielding (substituting $\partial_t \to \partial_x$ and integrating the conservative $(\tau,u,E)$ equations)
\be\label{zndprof_form}
\bar u=
1-\bar \tau,
\quad
\bar e= 
\frac{\bar \tau(\Gamma e_++1-\bar \tau)}{\Gamma},
\ee
$$
\begin{aligned}
\bar \tau&= \frac{ (\Gamma+1) (\Gamma e_++1)} {\Gamma +2} 
	- \frac{ \sqrt{ (\Gamma+1)^2 (\Gamma e_++1)^2
- \Gamma (\Gamma +2) ( 1+2(\Gamma +1)e_+ -2q(z-1))  } } {\Gamma +2}.
\end{aligned}
$$
The $\bar z$ component can then be solved via
$\bar z'=k\phi(c\bar e(\bar z))z  $ on $x<0$ (reaction zone).
A nonreactive ``Neumann shock'' at $x=0$ connects the ignited state at $x=0^-$ to a quiescent state at $x=0^+$
(for both of which $z=1$), and the profile remains constant thereafter, i.e., for all $x\geq 0^+$.
This corresponds to the physical picture of a gas-dynamical shock moving into an unburned, quiescent gas at 
$x\to +\infty$, which, its temperature being raised by compression of the shock, ignites and burns steadily, leaving
a ``reaction spike'' in its wake, with completely burned gas at $x\to -\infty$.

	\subsubsection{Viscous (rNS) profiles}

\begin{figure}[ht]
\centering 
\begin{tabular}{cc}
(a) \includegraphics[width=0.4\textwidth]{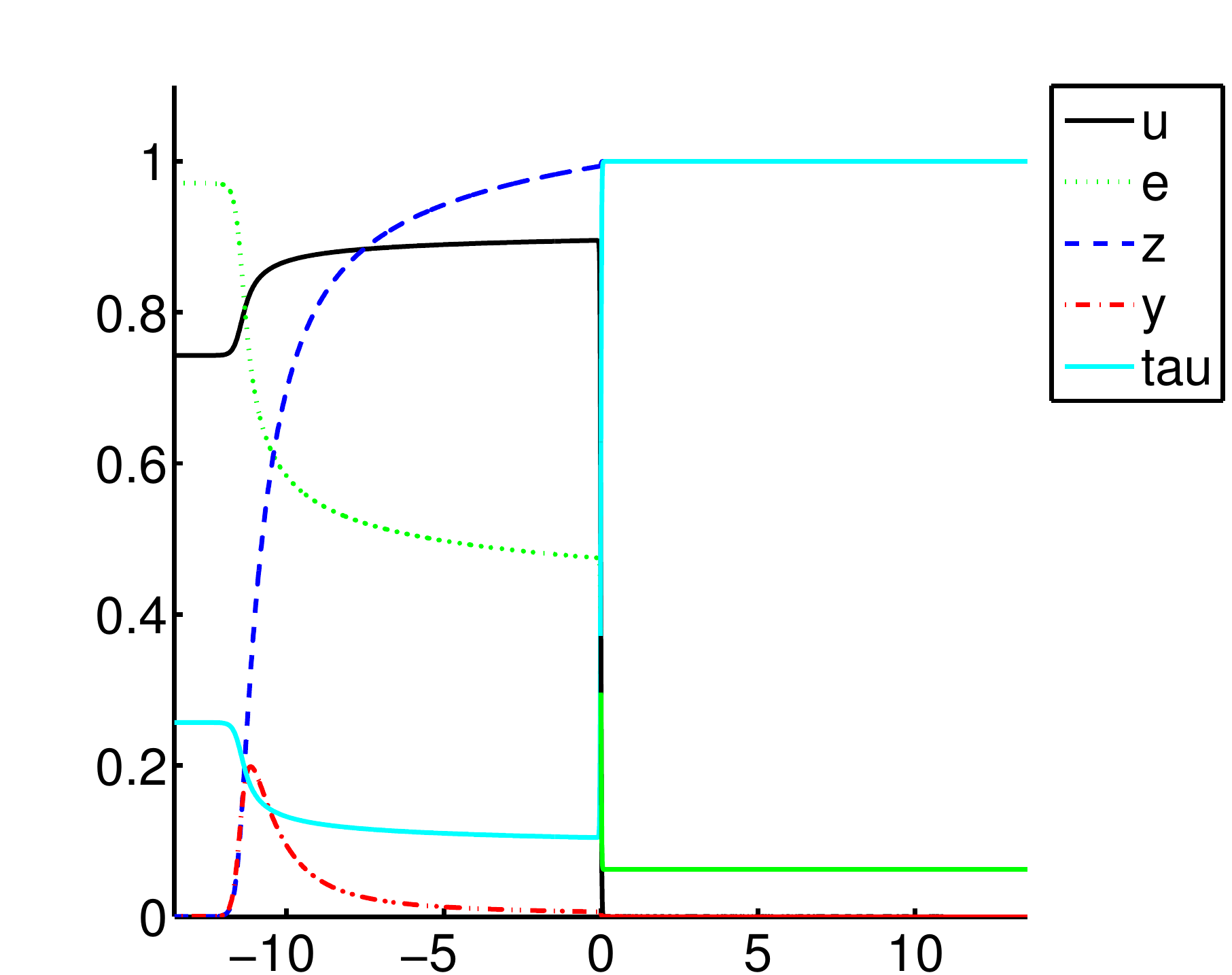} &
(b) \includegraphics[width=0.4\textwidth]{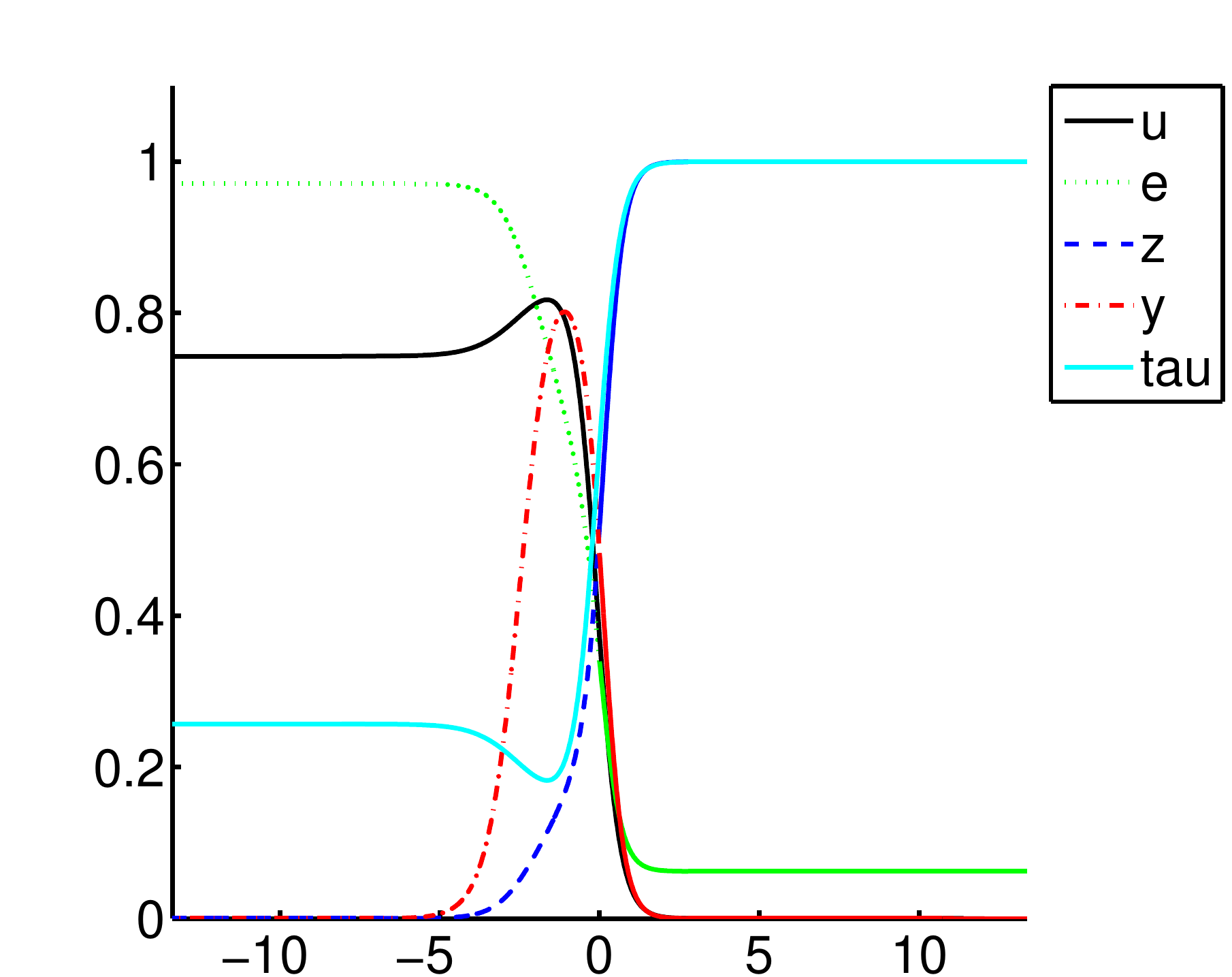}
\end{tabular}
\caption{Sample profiles illustrating diffusive effects. (a) $\nu=d=\kappa=0.01$. (b) $\nu=d=\kappa=0.3$.
In both cases the reaction zone structure is clearly visible, but in (b) the shock width is
of a similar order as the reaction zone width.  
%
For both plots, 
$e_\sp =6.23$e-2, $k = 1.53$e4, $q =6.23$e-1, $\Ea = 6$,
$\Gamma = 0.2$, $c_v= 1$. }
\label{fig:profiles}
\end{figure}

Likewise, parametrized by $(e_+,q, \mathcal{E}, \Gamma, \nu, \kappa, d) \in$ compact
domain (i.e., with nonphysical value $e_+=0$ adjoined), rNS profiles are exponentially convergent to their
endstates except at the degenerate ``Chapman--Jouget'' value $q=q_{CJ}$ \cite{LyZ1,Z2,Z3}, for which they decay algebraically.
Existence of rNS profiles for small viscosity/heat conduction/species diffusion
has been shown, for example, in \cite{GS,Wi}, by singular perturbation of the ZND case.
When diffusion coefficients are not small, profiles must be found in general numerically \cite{BHLyZ}.
Numerically determined profiles for different values of diffusion coefficients are displayed in Figure \ref{fig:profiles}.

\subsubsection{Issues and objectives}\label{s:stabprop}
\begin{figure}
\includegraphics[width=1.0\textwidth]{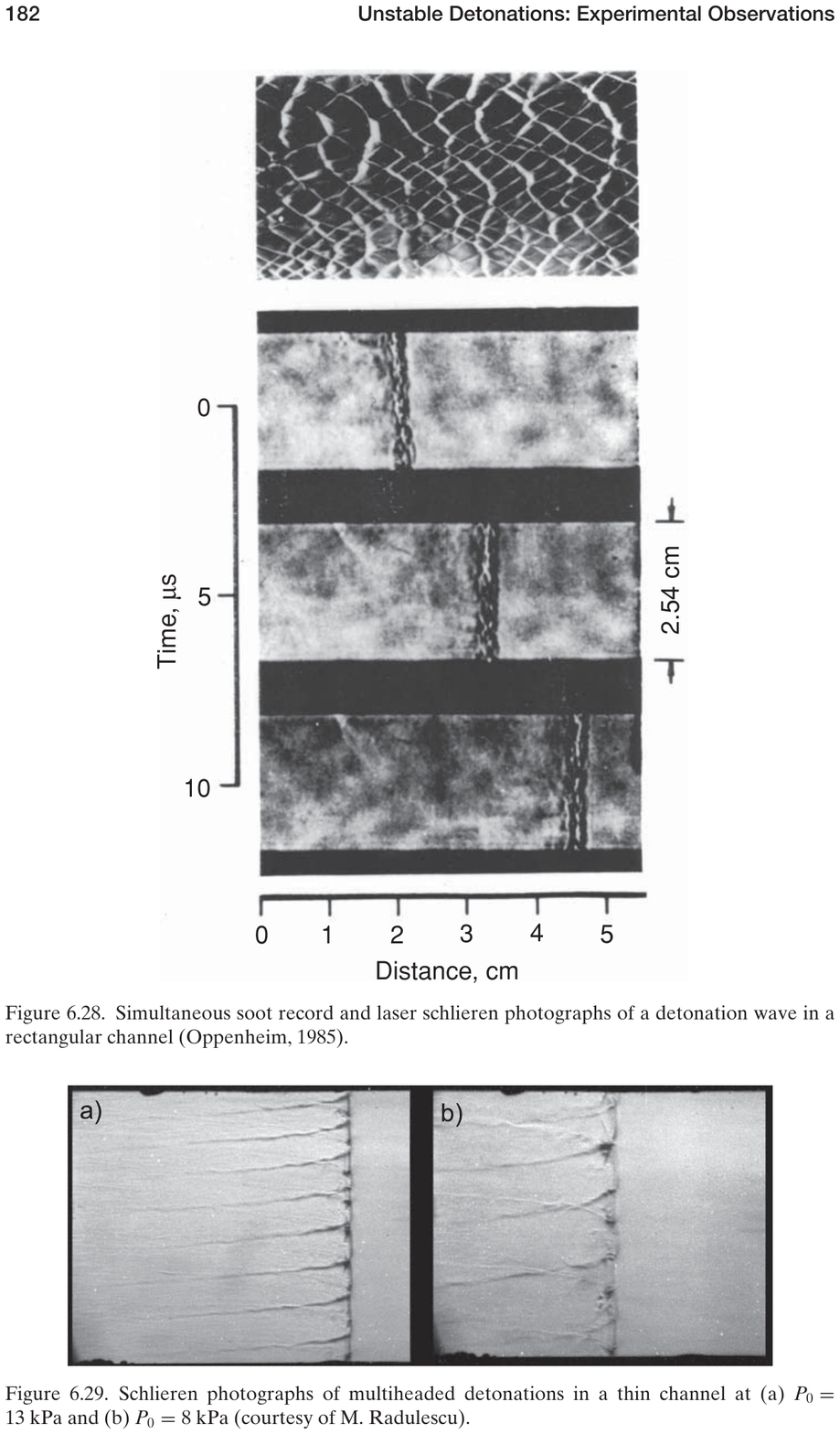}
\caption{ Detonation instability in a duct (John H.S
Lee, McGill University).}
\label{fig:Lee}
\end{figure}

Unlike nonreactive shocks, which are typically quite stable, detonations frequently exhibit instabilities of different types.
See Figure \ref{fig:Lee} depicting results of shock tube experiments carried out by John H.S. Lee (reprinted from \cite{L}
with permission of the author), which indicates the variety of possible behaviors as physical parameters are varied, from
a nonreactive shock-like coherent planar detonation layer, to apparent bifurcation to cellular or pulsating patterns,
to what appears to be chaotic flow.

The first mathematical model of detonation, the {\it Chapman--Jouget} model ($\sim$ 1890's; e.g., \cite{C,J1,J2})
treated detonations as a shock modified by instantaneous reaction.
This is sufficient to predict possible endstates and speeds
of planar discontinuities, but not to determine realizility by a connecting longitudinal reaction/dissipation structure.
Moreover, it does not capture the complicated instability/bifurcation phenomena described above; indeed, for the one-step
polytropic model of Example \ref{rNSeg}, Chapman-Jouget detonations are {\it universally stable} \cite{MaR,JLW}.

The modern theory of detonation stability dates from the post-world war II period, 
with the introduction of the ZND mdel \cite{Ze,vN1,vN2,D} and the pioneering stability/behavior studies of 
J.J. Erpenbeck and others.
The ZND model has successfully modeled a wide range of experimentally observed phenomena in stability/behavior.
Indeed, there is by now a comparatively long history ($\sim$ 1960's; e.g., \cite{Er1}), and extensive numerical
and analytical literature in the context of ZND; see, for example, \cite{Er1,LS,KS,CF,FD,BMR,B}, and references therein.
By contrast, until recently ($\sim$ 1990's; e.g. \cite{LyZ1}), there 
was relatively little investigation of the more complicated rNS model.

\medskip
{\bf Issues:} 
1. Experimental stability transitions/bifurcation to time-periodic pulsating/cellular wave patterns are 
well modeled by ZND.  But, there is no corresponding nonlinear stability or bifurcation theory,
and little regularity (or even well-posedness
for the (hyperbolic) equations.
2. The rNS equations on the other hand feature better regularity/well-posedness, 
but are significantly more complicated; till recently, 
there was neither linear data nor nonlinear theory.  Practical effects/importance of added transport
(viscosity/heat conduction/diffusion) terms is not clear.

\medskip

{\bf Objectives:} 
1. Review and rigorous (analytical) verification of conclusions plus systematic (numerical/analytical) 
exploration of parameter space; justification (and improvement) of numerics, for both (ZND) and (rNS). 
2. Systematic comparison between and synthesis of (rNS) and (ZND).

\subsection{Stability framework: normal modes analysis for ZND}\label{I.II}
(Following \cite{JLW}) Shifting to coordinates $\tilde x=x-st$ moving with the background
Neumann shock, write \eqref{rNS} as
$$
\begin{aligned}
W_t + F(W)_x=R(W),
\end{aligned}
$$
where
$$
\begin{aligned}
W:=\begin{pmatrix} u\\z \end{pmatrix}, \quad
F:=\begin{pmatrix} f(u)-su\\-sz \end{pmatrix},\quad
R:=\begin{pmatrix} qkz\phi(u)\\-kz\phi(u) \end{pmatrix}.\\
\end{aligned}
$$

\subsubsection{Fixed-boundary formulation}
Defining the Neumann shock location as $X(t)$, we reduce to a fixed-boundary
problem by the change of variables $x\to x-X(t)$.
In the new coordinates,
$$
\begin{aligned}
W_t + (F(W) - X'(t) W)_x=R(W), \quad x\ne 0,
\end{aligned}
$$
with jump condition
$$
\begin{aligned}
X'(t)[W] - [F(W)]=0, 
\end{aligned}
$$
$[h(x,t)]:=h(0^+,t)-h(0^-,t)$ denoting the jump at $x=0$ of a function $h$.

\subsubsection{Linearized equations}
Linearizing about
$(\bar W^0,0)$, we obtain the {\it linearized equations} 
$$
\begin{aligned}
(W_t - X'(t)(\bar W^0)'(x)) + (AW)_x =EW,
\end{aligned}
$$
$$
\begin{aligned}
X'(t)[\bar W^0] - [A W]=0, \quad x=0,
\end{aligned}
$$
where
$ A:= (\partial/\partial W)F$, $E:= (\partial/\partial W)R$. 

\subsubsection{Reduction to homogeneous form}

To eliminate the front from the interior equation, 
reverse the original transformation to linear order
by the change of dependent variables
$
W\to W- X(t)(\bar W^0)'(x),
$
motivated by $W(x-X(t),t))-W(x,t)\sim -X(t)W_x(x,t)
 \sim -X(t) (\bar W^0)'(x)$
 approximating to linear order the original, nonlinear transformation. (The trick of the ``:ood unknown'' of Alinhac \cite{A,JLW}.)
Substituting using 
$(A(\bar W^0)(\bar W^0)'(x))_x=E(\bar W^0)(\bar W^0)'(x)$
gives
\begin{equation}\label{flin}
\begin{aligned}
W_t + (AW)_x =EW
\end{aligned}
\end{equation}
with modified jump condition
$
X'(t)[\bar W^0]-  [A \big(W+ X(t) (\bar W^0)'\big) ]=0. 
$

\subsubsection{Generalized eigenvalue equation}

Seeking normal mode solutions $W(x,t)=e^{\lambda t}W(x)$,
$X(t)=e^{\lambda t}X$ yields the generalized eigenvalue equations
$
(AW)' = (-\lambda I   + E)W, \quad x\ne 0,
$
$
 X(\lambda[\bar W^0]-[A (\bar W^0)']) - [A W]=0, 
$
where ``$\prime$'' denotes $d/dx$. or, setting $Z:=AW$, to
\be\label{int}
Z' = GZ, \quad x\ne 0,
\ee
\ba\label{bound}
X(\lambda[\bar W^0]-[A (\bar W^0)']) - [Z]=0, \quad x=0, 
\ea
where $ G:=(-\lambda I  + E)A^{-1}. $

\subsubsection{Stability determinant}

We define the {\it Evans--Lopatinski determinant}
\begin{equation}
\begin{aligned}\label{fform}
D_{ZND}(\lambda)&:=
\det \begin{pmatrix}
Z_1^-(\lambda,0), & \cdots, & Z_{n}^-(\lambda,0), &
\lambda [\bar W^0]-[A (\bar W^0)']\\
\end{pmatrix}\\
&=
\det \begin{pmatrix}
Z_1^-(\lambda,0), & \cdots, & Z_{n}^-(\lambda,0), &
\lambda[\bar W^0]+A(\bar W^0)'(0^-)\\
\end{pmatrix},\\
\end{aligned}
\end{equation}
where $Z^-_j(\lambda,x)$ are a basis of solutions of the interior equations \eqref{int} 
decaying as $ x\to -\infty$.
By $ A(\bar W^0)':=dF(\bar W^0)(\bar W^0)'=R(\bar W^0)$ plus duality, we can rewrite \eqref{fform} in the simpler form
$$
D_{ZND}(\lambda)= \tilde Z_n^-(\lambda,0)\cdot
\big(\lambda[\bar W^0]+R(\bar W^0)(0^-)\Big)
$$
useful for numerics \cite{Br,HuZ} and also analysis \cite{Z1,Z2},
where $\tilde Z_n^-$ is a (unique up to constant multiple) solution of the dual equation
$\tilde Z'=-G^* \tilde Z$ decaying as $x\to -\infty$.
The function $D_{ZND}$ is exactly the {\it stability function} derived in a different form by Erpenbeck \cite{Er1,BZ1}.

\medskip

$\bullet$ {\it Evidently, $\lambda$ is a generalized eigenvalue iff $D_{ZND}(\lambda)=0$.}

\begin{definition}\label{zndss}
A ZND profile is spectrally stable if there are no zeros of the associoated Lopatinski determinant in 
$\{\lambda:\; \Re \lambda \geq 0\} \setminus \{0\}$ \cite{Er1}.
(By translation-invariance, there is always a zero at $\lambda=0$.)
\end{definition}

\subsection{Normal modes analysis for rNS}

Take without loss of generality $s=0$ (co-moving coordinates), $\nu=1$, so that
$u=\bar u(x)$ is an equilibrium.
Linearized eigenvalue equations
$$
\lambda W=LW:= -(A(x)W)_x +\eps (B(x) W_x)_x +EW
$$
may be written as a first-order system
$$
\mathcal{Z}'=\mathcal{A}(x,\lambda)\mathcal{Z},
$$
where $\mathcal{Z}=
\begin{pmatrix} Y\\ W_2 \end{pmatrix}=
\begin{pmatrix} AW-\eps BW_x\\ W_2 \end{pmatrix} $
	is an augmented ``flux'' variable \cite{Z2}.

Define the {\it Evans function}
\be\label{Ddef1}
D_{rNS}(\lambda):=\det (\mathcal{Z}_1^-, \dots, \mathcal{Z}_k^-, 
\mathcal{Z}_{k+1}^+, \dots, \mathcal{W_N})|_{x=0}
\ee
where $ \{\mathcal{Z}_1^-, \dots, \mathcal{Z}_k^-\}(\lambda,x)$ 
and $\{ \mathcal{Z}_{k+1}^+, \dots, \mathcal{Z}_N)\}(\lambda,x)$ are bases of
solutions decaying as $x\to \infty$ and $x\to +\infty$.

\medskip
$\bullet$ {\it Evidently, $\lambda$ is an eigenvalue iff $D_{rNS}(\lambda)=0$.}

\begin{definition}\label{rnsss}
An rNS profile is spectrally stable if there are no zeros of the associoated Evans function in
$\{\lambda:\; \Re \lambda \geq 0\} \setminus \{0\}$ \cite{LyZ1,LRTZ,TZ4}. (As for ZND, 
there is always a zero at $\lambda=0$.)
\end{definition}

\subsection{Abstract viscous stability results}\label{1.III}
Let $\{\bar W^\e\}$ be a one-parameter family of 
viscous strong detonation waves for rNS with polytropic equation of state \eqref{eoseg}.

\subsubsection{Spectral stability transitions}\label{transitions}

\begin{lemma}[Stability in the small-heat release limit \cite{LyZ1}]\label{qlem}
If $q^\e\to 0$ as $\eps\to 0$, then 
$\{\bar W^\e\}$ is spectrally stable for $\eps$ sufficiently small.
\end{lemma}

\begin{lemma}[Absence of steady bifurcations \cite{LyZ1}]\label{biflem}
	For all $\e,$ the associated Evans function has a zero of multiplicity
one at $\lambda=0:$ $ D(\e,0)=0$, and $D'(\e,0)\ne 0$, hence stability transitions if they occur involve passage of
nonzero conjugate zeros across the imaginary axis.
More generally, this holds for any equation of state for which the associated CJ profiles are stable

\end{lemma}

Lemma \ref{qlem} is an immediate consequence of the construction of an Evans function, done similarly as in \cite{GZ,AGJ},
using the resulting continuity with respect to parameters together with decoupling at $q=0$ of gas-dynamical
($u$) and reaction ($z$) equations.
Lemma \ref{biflem} follows by a ``stability index'' computation like those of \cite{PeW,GZ,ZS}, quantifying the
intuition that low-frequency behavior of rNS should ``not see'' reaction and transport scales, so shoul
reduce to that of CJ.
See \cite{JLW} for a far-reaching extension of this principle including also ZND and multi-D.

\medskip

{\it Consequences:} 1. Stable waves exist. 2. Stability transitions should they occur
are of (spectral) Hopf, i.e., ``pulsating'' type, {\it as seen in experiment.}
(Link between behavior and equations.)

\subsubsection{Nonlinear stability/bifurcation criteria}\label{s:crit}

\begin{theorem}[Spectral $\Rightarrow$ nonlinear stability \cite{TZ4}]\label{stabthm}
For all $\e,$ $\bar W^\e$ 
is $L^1\cap L^p\to L^p$ linearly orbitally stable if and only if, for all $\e,$
the only zero of 
$D(\e,\cdot)$ in $\Re \lambda \ge 0$ is a simple zero at the origin,
in which case $\bar W^\e$ is $L^1 \cap H^3 \to L^1 \cap H^3$ linearly 
and nonlinearly orbitally stable, with 
$$
|\tilde W^\e(\cdot, t)-\bar W^\e(\cdot -\delta(t))|_{L^p}
\le C|\tilde W^\e_0-\bar W^\e|_{L^1\cap H^3}(1+t)^{- \frac{1}{2}(1-\frac{1}{p})},
$$
for nearby solutions $\tilde W^\e$, where 
$$
\begin{aligned}
|\delta(t)| &\le C|\tilde W^\e_0-\bar W^\e|_{L^1\cap H^3},
	\quad
|\dot \delta(t)| \le C|\tilde W^\e_0-\bar W^\e|_{L^1\cap H^3}(1+t)^{-\frac{1}{2}}.
\end{aligned}
$$
\end{theorem}

\begin{theorem}[Spectral $\Rightarrow$ nonlinear bifurcation \cite{TZ4}]\label{bifthm}
Assume that $\bar W^\e$ undergoes transition from linear stability to linear instability at $\e = 0,$
via passage of a single complex conjugate pair of eigenvalues 
$\lambda_\pm (\e)=\gamma(\e)+i\tau(\e)$
through the imaginary axis:
\be\label{nondeg}
\gamma(0)=0, \quad \tau(0)\ne0, \quad d\gamma/d\e(0)\ne 0.
\ee
Then, given exponential weight $\o$, for $0\leq r \ll 1$ and $C\gg 1$, 
there are $C^1$ functions $r \to \e(r), T(r),$ with $\e(0)=0$, $T(0)=2\pi/\tau(0),$
and a $C^1$ family of time-periodic solutions  $\tilde U^{r}(x,t)$
of (rNS) with $\e=\e(r)$, of period $T(r)$,
with
$
C^{-1}r \le \|\tilde U^r-\bar U^\eps\|_{H^2_\o}\le Cr.
$
Up to translation in $x$, $t$,
these are locally unique in $\|\cdot\|_{H^2_\o}$.
\end{theorem}

Theorem \ref{stabthm} is established by detailed pointwise Green bounds obtainedf from stationary phase type estimates
on the inverse Laplace transform representation of the linearized solution operator, together with
a nonlinear shock tracking argument, in the spirit of \cite{ZH,MaZ,Z4}.
Theorem \ref{bifthm} is established by a novel ``reverse temporal dynamics'' argument 
using inverse Laplace transform estimates similar to those for stability.  See also \cite{TZ2,TZ3,SS} for related
studies in the shock wave case.  For a nonlinear stability analysis of the bifurcating time-periodic solutions,
see \cite{BSZ}.

\medskip

{\it Consequences:} Spectral information as in Lemmas \ref{qlem}--\ref{biflem} translates to full nonlinear results.

\subsubsection{Closing the philosophical loop: the rNS$\to$ZND limit}\label{s:loop}
At this point, the situation as regards the two theories (rNS and ZND) is that we have for ZND
decades of spectral stability data, numerics, and formal asymptotics for ZND, but no nonlinear theory; for rNS,
we have essentially the reverse.
A way to repair this situation, combining the strengths of the two theories, is to link them via the vanishing viscosity,
rNS$\to$ZND limit.
The limiting profile structure problem has been studied in \cite{GS,Wi}, etc., with definitive results.
However, until recently, the only analytical result regarding stability was the study 
by Roquejoffre--Vila \cite{RV} for {\it Majda's model} \cite{Ma1}, a simplified qualitative model of detonations. 
A generalization to the full rNS system is as follows; here, $\bar W^\eps=$ represents an $\eps$-profile,
with $\eps$ measuring size of transport (viscosity/heat conduction/diffusion) coefficients

\begin{theorem}[rNS spectrum in the ZND limit \cite{Z2}]\label{limthm}
Spectral stability of $\bar W^\eps$
for $\eps>0$ sufficiently small
is equivalent to spectral stability of the limiting ZND 
detonation $\bar W^0$ together
with spectral stability of the viscous version of the associated Neumann shock.
Moreover, 
(i) For $C\le |\lambda| \le C/\eps$, $C$ sufficiently large,
on $\Re \lambda>-\eta$ for $\eta$, $\eps>0$ sufficiently small, 
$\eps$ times each zero of $D^\eps_{rNS}$ converges
to a zero of $\frac{D_{NS}(\lambda)}{\lambda}$ on $\Re \lambda \ge 0$
and each zero of $D_{NS}$ on $\Re \lambda>0$ is the limit of $\eps$
times a zero of $D_{rNS}$  on $\Re \lambda>0$, for $C\le |\lambda| \le C/\eps$.
(ii) For $|\lambda|\le C_0$, $C_0$ arbitrary, on $\Re \lambda \ge -\eta<0$,
the zeros of $D^\eps_{rNS}$ converge
in location/multiplicity as $\eps\to 0$ to the zeros of $D_{ZND}$.
\end{theorem}

The proof of Theorem \ref{limthm} is by detailed multi-scale analysis as in stability of strong shocks
and other asymptotic limits \cite{PZ,HLZ}, together with an $\eps$-variational argument like that used 
in \cite{GZ} and \cite{ZS} to study the related low-frequency (small-$\lambda$) limit.
The detailed asymptotics provided on the profile by the analyses of \cite{GS,Wi} are used in an important way.
It is known that nonreactive viscous shocks of a polytropic gas are universally stable \cite{HLyZ1,HLyZ2},
hence the theorem {\it reduces spectral stability for rNS in the small-viscosity limit to spectral stability of ZND}.

\medskip

{\it Consequences:} 
1. Verifies NS stability/bifurcation for small $\eps$ through extensive existing numerical studies for ZND.
2. Gives rigorous nonlinear sense to (spectral) ZND results.

\medskip

This gives one answer to the question ``what is the
role of viscosity?'' (namely, logical development/foundations).
We'll explore a possible different answer below, in Section \ref{I.VI}.

\subsection{Abstract inviscid stability results}\label{1.IV}
(First rigorous stability results for ZND)
Let $\{\bar W^\e\}$ be a one-parameter family of 
strong detonation waves for ZND with polytropic equation of state \eqref{eoseg}.

To explain our next results, we first recall that the parametrization given in Section \ref{s:param}
is not the standard one given in the literature, but our own ``improved'' version \cite{Z2}.
In the classical parametrization given e.g. in \cite{Er2}, $e_+$ rather than speed $s$ is held fixed,
and the detonation parametrized rather
by the {\it overdrive} $1<f<\infty$, defined as
the square of the ratio of relative speed of the detonation (with
respect to the ambient gas) and the minimum, Chapman--Jouget, detonation
speed among all possible strong detonations \cite{Er2,FW,LS,BMR}.
In this classical scaling, two
rules of thumb observed numerically are that detonations are more stable the smaller the heat release $q$
and the higher the overdrive $f$.
The former was proved by Erpenbeck for finite frequencies, but his treatment of high frequencies was incomplete \cite{Z1}.

\begin{lemma}[Stability in the small-heat release limit \cite{Z1}]\label{qzndlem}
	In the scaling of Section \ref{s:param},
if $q^\e\to 0$ as $\eps\to 0$, then $\{\bar W^\e\}$ is spectrally stable for $\eps$ sufficiently small.
\end{lemma}

\bc [Stability in the high-overdrive limit \cite{Z1}]\label{odcor}
In the scaling of Erpenbeck \cite{Er2}, ZND detonations of \eqref{eoseg} are 
spectrally stable in the fixed-activation energy, fixed-heat release, high-overdrive limit $f\to \infty$.
\ec

The first result includes but is not restricted to the observation of Erpenbeck that, in the scaling of \cite{Er2},
ZND detonations are stable in the fixed-activation energy, fixed-overdrive, small-heat release limit,
which in our scaling corresponds to fixed-activation energy, fixed-$e_+$ or shock strength, and small=$q$ or
heat release.
The second result, corresponding in our scaling to stability in 
the simultaneous {\it zero-heat release}, {\it zero-activation energy $\mathcal{E}$},
and {\it strong-shock} (zero-$e_+$) {\it limits},
resolves an open problem from \cite{Er2}.
Our favorable coordinatization ($s=1$ held fixed, $e_+\to 0$), suggested by similar scalings used to study the strong-shock
limit for gas dynamics \cite{HLZ,HLyZ1}, plays an important role in the analysis.  For, this keeps all quantities
bounded for bounded frequencies, independent of parameters, allowing uniform treatment of the strong-shock limit.
By contrast, internal energy $e$ and temperature $T$ blow up for the classical scaling in the strong-shock limit.
Results obtained in passing are 1D high-frequency stability, new asymptotic ODE techniques.

\medskip

{\it Consequences:} 
1. Analytical signposts guiding delicate/computationally intensive numerics \cite{Er1,LS}.
2. 1D high-frequency stability, validating numerics by truncation of computational domain.

\br\label{hf1rmk}
The 1-D high-frequency stability analysis foreshadows issues addressed in Section \ref{II.I}
for multi-D.  {\it Notably, the 1D analysis requires only $C^2$ regularity on coefficients/equation of state.}
\er

\subsection{Numerical results for ZND}\label{I.V}

\subsubsection{Natural coordinatization}\label{nat}
\begin{figure}
\includegraphics[trim=0 30 0 30, scale=.45]{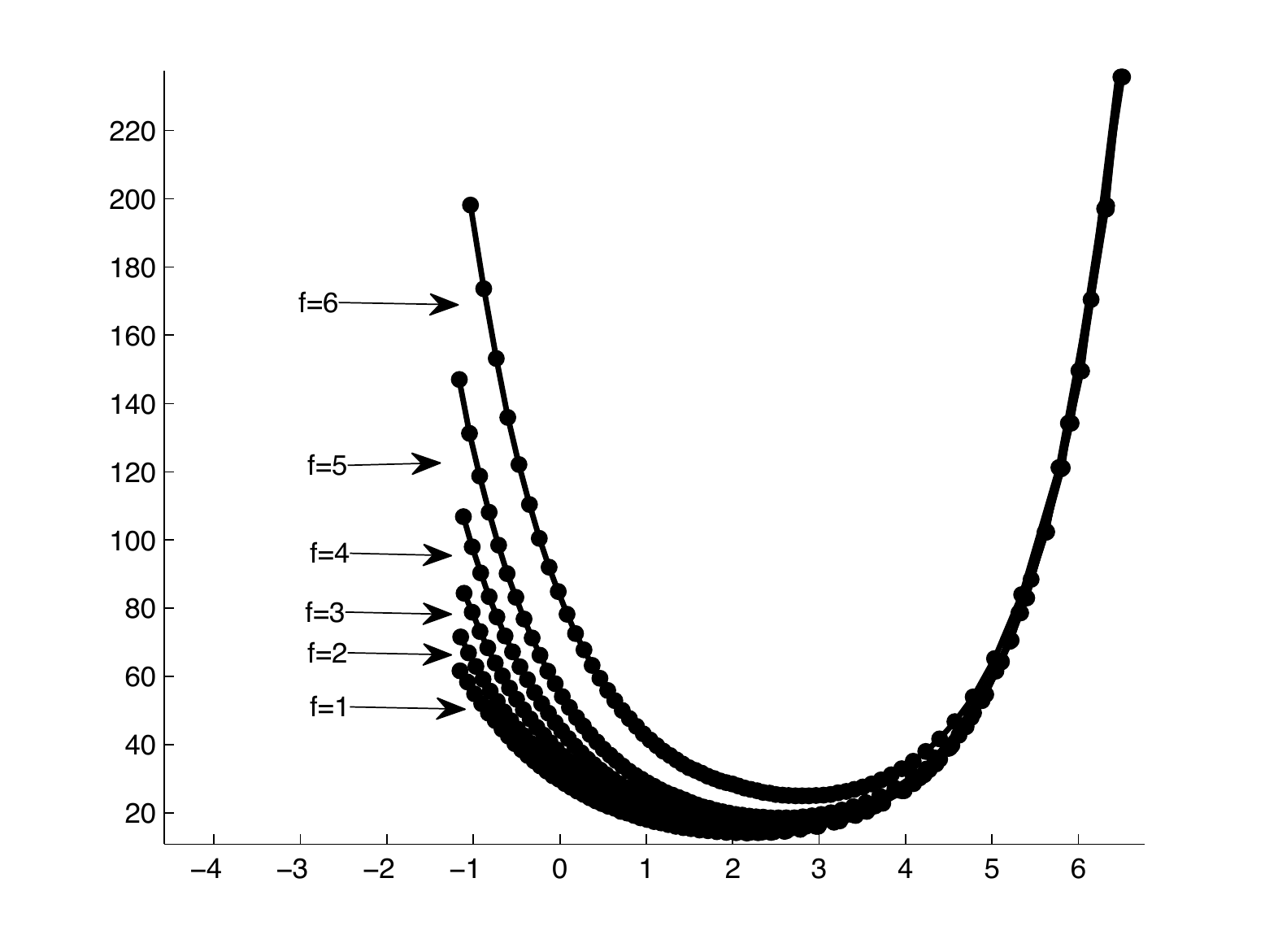}
\caption{ $\Gamma=0.2$, constant overdrive $f$, $\E_0$ vs. $\log Q_0$.}
\label{fig:overdrive}
\end{figure}

The novel scaling introduced in Section \ref{s:param} is helpful not only 
for rigorous analysis, as seen in Section \ref{1.IV}, but also at the level of numerics/modeling.
In Figure \ref{fig:overdrive}, we display in the classical scaling of Erpenbeck \cite{Er2}
results for a standard benchmark problem
of Fickett and Woods \cite{Er2,FW,LS}, holding overdrive $f$ fixed and varying activation
energy $\mathcal{E}_0$ and heat release $Q_0$, with $\Gamma=1.2$.
The solid curves depicted are the {\it neutral stability curves} across which detonations change from stable (below) 
to unstable (above) as $\mathcal{E}_0$ is increased.
In this figure, we see the stabilizing effect of increasing $f$ and the destabilizing effect of increasing $\mathcal{E}_0$;
however, there is an apparent hysteresis effect as $\mathcal{Q_0}$ is increased, with detonations first destabilizing, then
restabilizing for large $\mathcal{Q_0}$.
Moreover, there is a singularity at the right of the diagram with $\mathcal{E}_0, Q_0\to \infty$.

\begin{figure}
\includegraphics[trim=0 200 0 200, scale=.45]{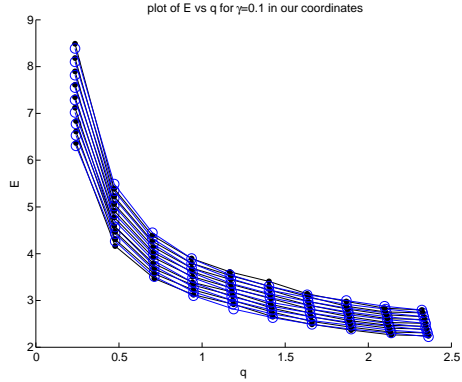}
\caption{$\cE=\E_0 e_+$ vs. $q=Q_0 e_+$; polynomial fit, average relative error $1\%$.}
\label{fig:our_coordinates}
\end{figure}

In Figure \ref{fig:our_coordinates}, we depict the analogous neutral stability curves for the same gas constant
$\Gamma=1.2$ in {\it our} scaling (the one of Section \ref{s:param}), 
with $e_+$ held fixed and $\cE=\E_0 e_+$ and $q=Q_0 e_+$ varying.
In these coordinates, both hysteresis and singularity are removed.
The latter allows us to verify numerically {\it stability in the zero-activation energy limit}:
$\mathcal{E}=0$ $\Rightarrow$ ZND stability for any $q$, $e_+$.

\begin{figure}
\includegraphics[trim=0 200 0 200, scale=.45]{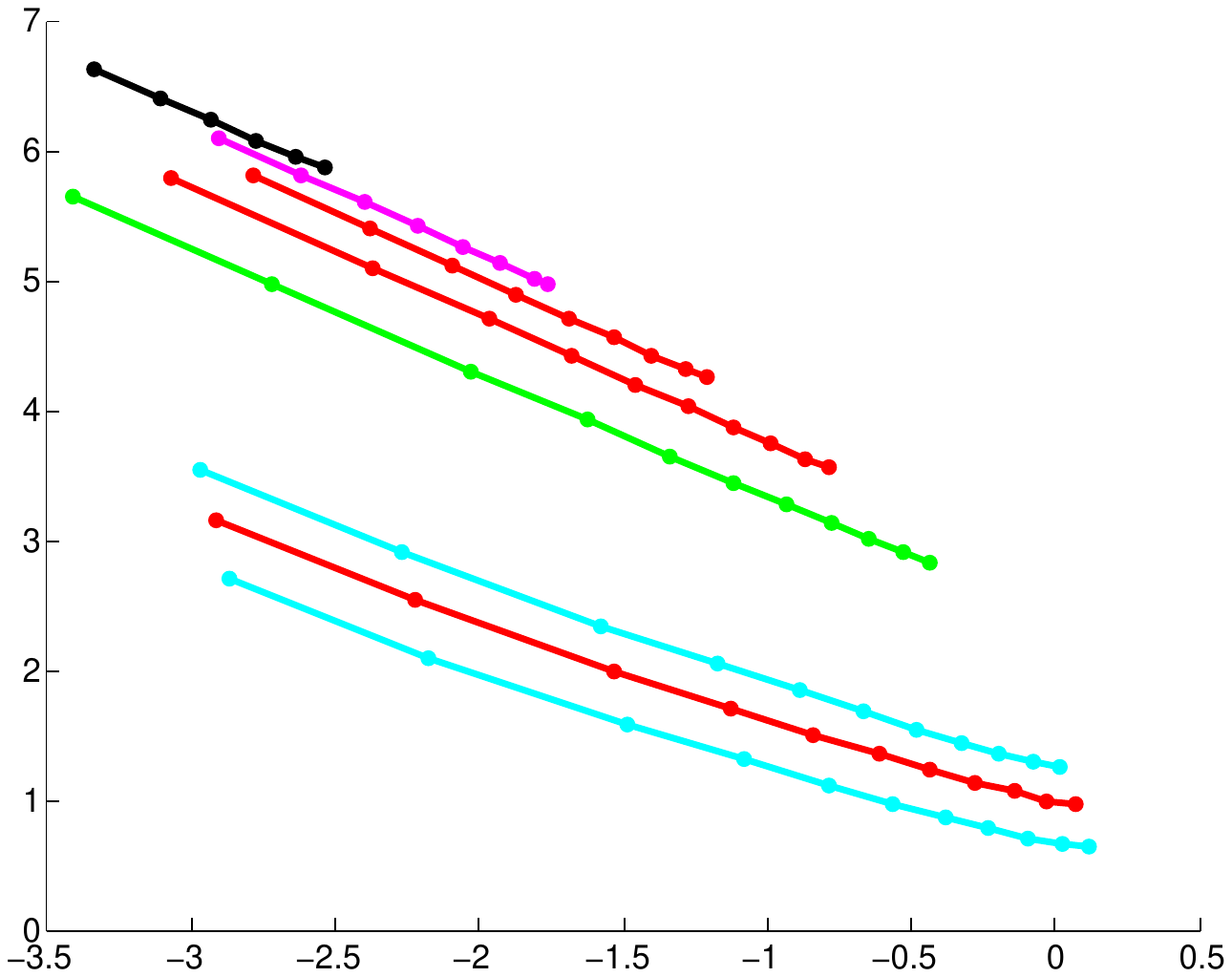}
\caption{ $\log \cE$ vs. $\log q$.}
\label{fig:fitcurves}
\end{figure}

Moreover, the neutral stability curves follow a simple and regular pattern, as may be seen most
dramatically in the log-log plot of Figure \ref{fig:fitcurves}.
Indeed, a naive polynomial fit with 20 stored coefficients 
is sufficient to recover the entire diagram in seconds with $2\%$ minimum/$1\%$ average accuracy,
a considerable compression of data for a diagram that required a reported
5 hours on a Cray supercomputer in 1990 to produce a single fixed-overdrive curve \cite{LS}.

\subsubsection{Computational improvements}\label{s:comp}
Besides the improvement in parametrization described above, we have 
by adapting to detonation theory numerical Evans function algorithms developed for the study of viscous shock stability \cite{BHZ},
improved computation speed by a factor of 1-2 orders of magnitude compared to the current state
of the art as described, e.g., in \cite{LS,SK}; see \cite{HuZ2,BZ1}.
With these improvements, combined with vastly improved hardware capability,
what took 5 hours on a supercomputer in 1990 to compute a single fixed-overdrive curve
today takes 5 hours on a Mac Quad Duo to compute the full Figure \ref{fig:overdrive}.
Indeed, this can be carried out perfectly well on a laptop.

We are now able to not only compute neutral stability curves, but to accurately describe all unstable eigenvalues
even for large activation energies; see for example the eigenvalue configuration displayed in
Figure \ref{fig:bench}(LEFT) for the same benchmark problem studied in Figures \ref{fig:overdrive}--\ref{fig:our_coordinates}
at activation energy $\approx 7.1$, for which we accurately resolve a pattern of $\approx 50$ unstable roots using
code supported in the MATLAB-based openware package STABLAB \cite{BHZ}.

\subsection{Numerical results for rNS}\label{I.VI}
Improvements in computions/power have made possible for the first time numerical Evans investigations for rNS,
a substantially more intensive problem than ZND.
These investigations, though just beginning, have already yielded surprising results

\subsubsection {``Viscous hyperstabilization''/hysteresis} \label{hyper} 
\begin{figure}[ht]
\centering
\begin{tabular}{cc}
\includegraphics[width=0.4\textwidth]{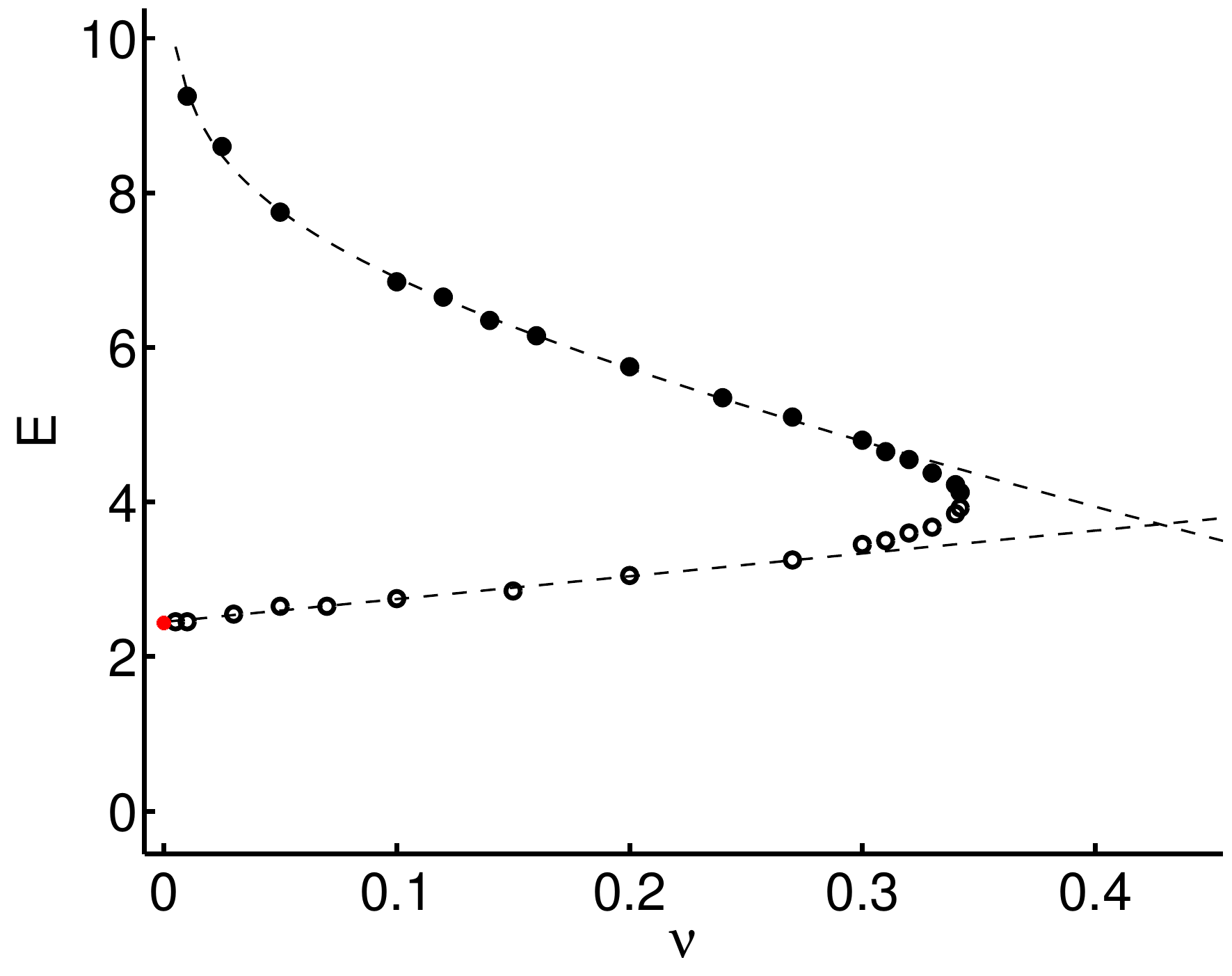}
\includegraphics[width=0.4\textwidth]{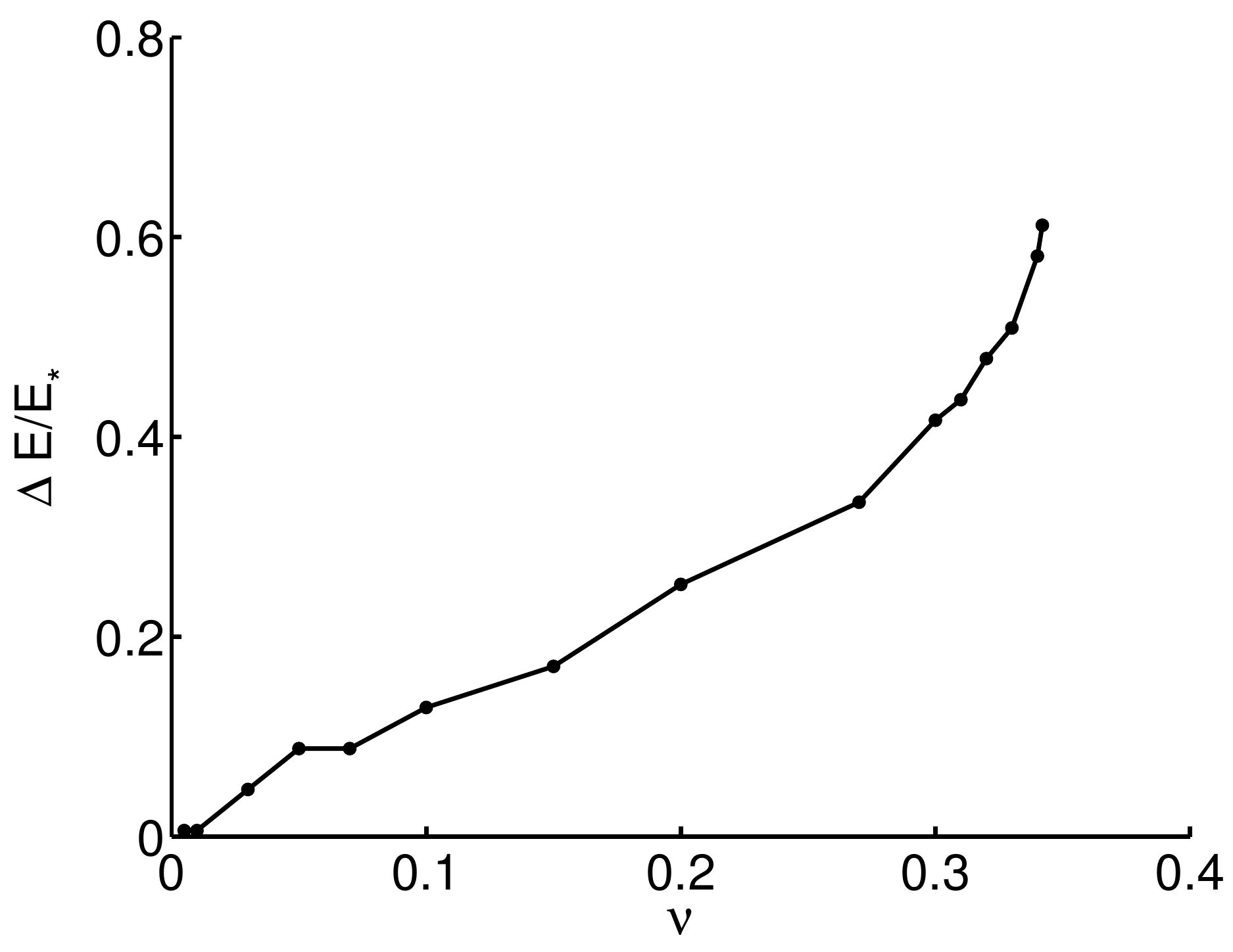}
\end{tabular}
\caption{LEFT: Neutral Stability Boundaries in the $\Ea$-$\nu$ plane.
	The best-fit curve (dashed line, $\nu<0.27$) for the upper boundary is $ \Ea^\sp(\nu)=5.67-6.16\nu -0.804 \ln(\nu) $. 
	RIGHT: Viscous delay (cf. \cite{RAP1,RAP2}): We plot $\Delta E/\mathcal{E}_*=(\Ea^\sm(\nu)-\mathcal{E}_*)/\mathcal{E}_*$ against $\nu$, where $\mathcal{E_*}$ is the approximation to the ZND neutral boundary. 
Here, $\nu  = d = \kappa$, $\Gamma = 0.2$, $e_\sp=6.23$e-2, and $q$ = 6.23e-1
are held fixed;
the red dot denotes the ZND (inviscid) stability boundary (lower boundary only!).
}
\label{fig:sb}
\end{figure}

For the benchmark problem discussed in Section \ref{I.V},
Romick et al \cite{RAP1,RAP2} have carried out numerical time-evolution studies indicating a significant
delay in transition to instability as activation energy is increased for the viscous (rNS) problem as compared
to the inviscid (ZND) one, as much as $10\%$ for values of viscosity in the high range of physically relevant scales.
Our numerical Evans investigations both confirm and extend these observations, indicating not only the expected delay
but also a new type of hysteresis in which viscous detonations 
{\it eventually restabilize} as activation energy is increased still further \cite{BHLyZ}.
This striking phenomenon is depicted in Figure \ref{fig:sb}(LEFT); see Figure \ref{fig:sb}(RIGHT) for a graph of
viscous delay vs. viscosity.
We call this phenomenon {\it viscous hyperstabilization}; we have conjectured \cite{BHLyZ} 
that it occurs for any nonzero viscosity, no matter how small.

Note the slow, apparently logarithmic, growth, in the upper stability boundary of Figure \ref{fig:sb}(LEFT) as viscosity goes
to zero, suggesting that hyperstablization might play a relevant physical role even for quite small values of viscosity.
Another notable feature of Figure \ref{fig:sb}(LEFT) is the ``nose'' to the right of the neutral stability curve,
where upper and lower boundaries meet. This indicates that there is {\it no instability}, regardless of the value of $\E$,
for sufficiently large viscosity.
For reference, the viscosity values considered in \cite{RAP1,RAP2} correspond to $\nu=0.1$ in the scaling of
Figure \ref{fig:sb}.

\subsubsection{Associated eigenvalue distributions}\label{s:dist}
The restabilization phenomenon just described is the more remarkable given the details of the unstable eigenvalue distribution.
In the inviscid case, it is more or less a universal principle that increasing $\E$ increases instability \cite{Er2,FD,LS}; 
indeed, as $\E$ increases, more and more unstable eigenvalues cross the imaginary axis from stable to unstable 
complex half-plane never to return, in a cascade of Hopf bifurcations.

\begin{figure}[ht] 
   \centering
   \begin{tabular}{cc}
   \includegraphics[width=.4\textwidth]{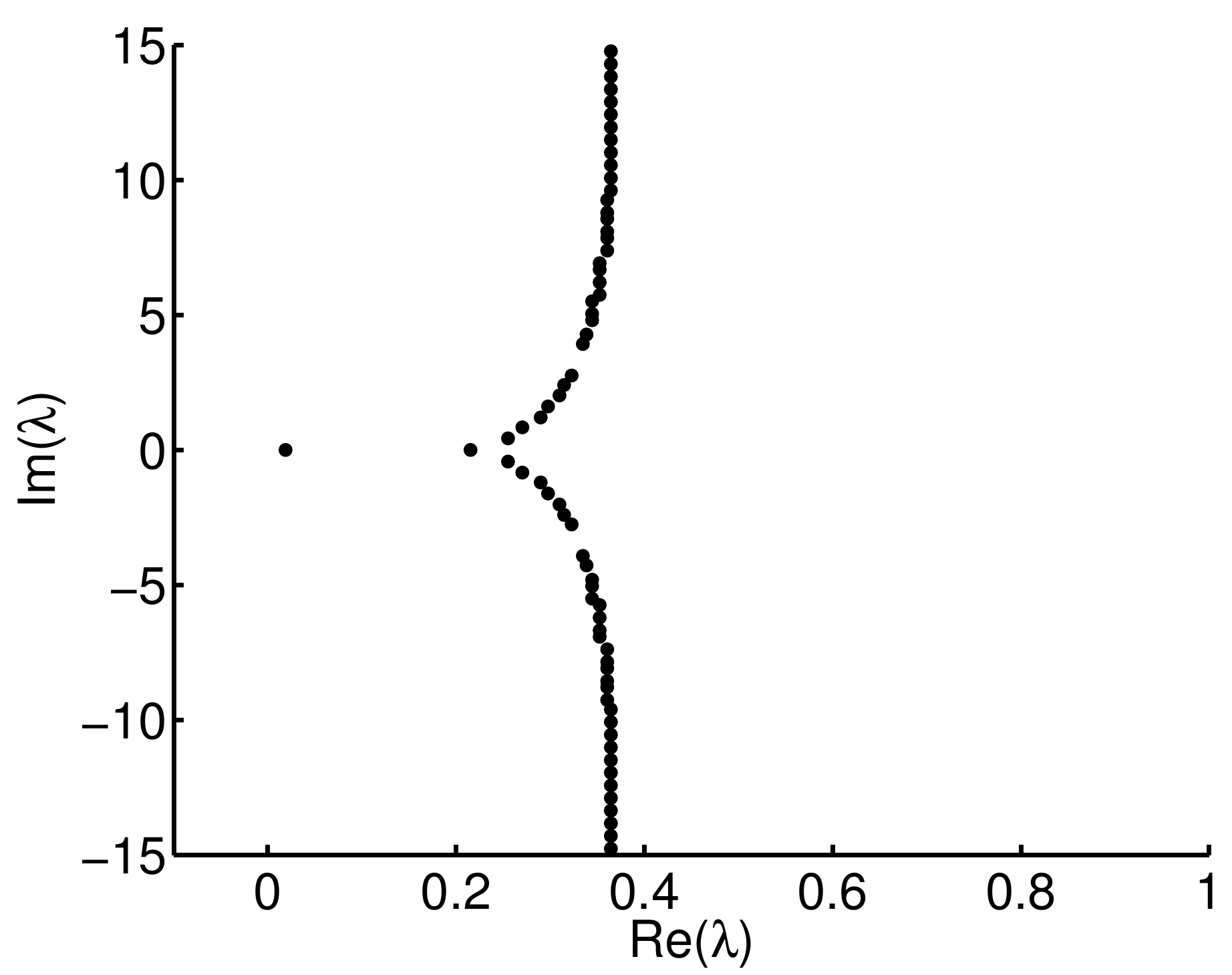} 
    \includegraphics[width=2.6in]{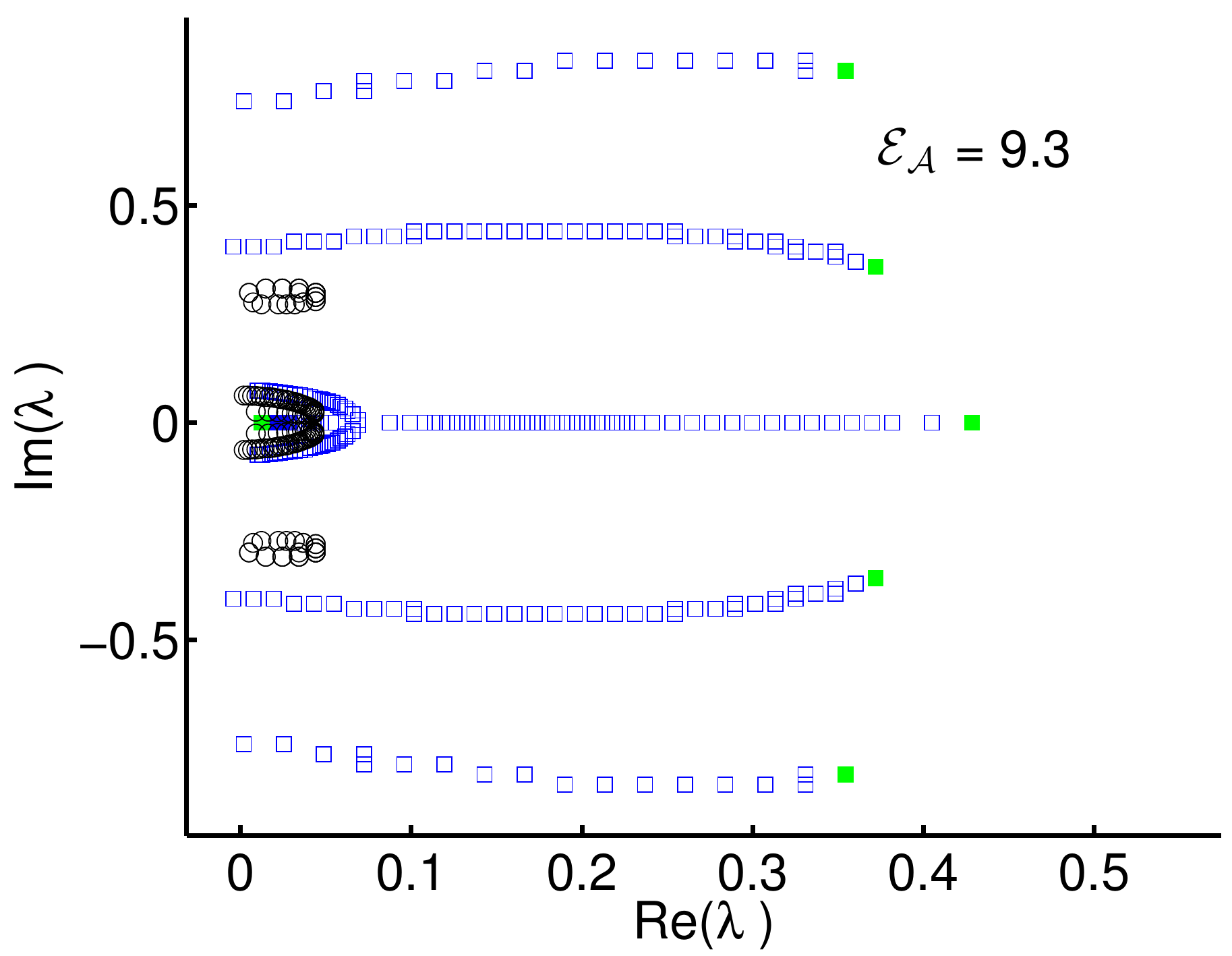} \\
   \end{tabular}
   \caption{LEFT: Unstable eigenvalues for ZND, $\Ea=7.1$; none for rNS!
	   RIGHT: 
    The movement of unstable roots in the complex plane as $\Ea$ increases. Circles 
    mark rNS roots, open squares 
    ZND roots.  The smaller modulus roots enters for $\Ea \approx 2.75$ (Panel (a)) and the higher modulus roots for $\Ea \approx 3.65$ (Panel (d)). The high modulus roots have a turning point at about $\Ea \approx 5.2$, and the smaller modulus roots have theirs at $\Ea \approx5.5$ (Panels (i) and (j)). The large modulus roots leave at $\Ea\approx 6.55$ and the small modulus roots leave at approximately $\Ea\approx6.85$.
   }
   \label{fig:bench}
\end{figure}

In Figure \ref{fig:bench}(LEFT) we display the eigenvalue distribution at $\E\approx 7.1$, for which there are 48 unstable
roots together with the translational eigenvalue at $\lambda=0$; further increases in $\E$ lead to further unstable
eigenvalues.
In Figure\ref{fig:bench}(RIGHT) we display for contrast the behavior of rNS eigenvalues for the value of viscosity
$\nu=0.1$ considered in \cite{RAP1,RAP2}, tracking the unstable eigenvalues as $\E$ is varied through the stability 
transition region.  For this viscous case, we find that there are {\it just two} pairs of unstable eigenvalues in total,
which after crossing the imaginary axis to the right turn back and rather quickly restabilize by crossing back into the 
stable half-plane; meanwhile, the nearby inviscid eigenvalues
plotted in the same figure may be seen to continue to the right.
At the value $\E\approx 7.1$ corresponding to the display of unstable inviscid eigenvalues in
Figure \ref{fig:bench}(LEFT), there are {\it no remaining} unstable eigenvalues for the viscous case with $\nu=0.1$.

\medskip

{\it Consequences:} 
1. Another possible answer to the question ``what is the role of viscosity'' in stability of strong detonations.
2. Potentially important physical effect, meriting further study.  

\subsection{Discussion and open problems}\label{disc1}
Investigation of detonation stability has proceeded by a blend of rigorous analysis, formal asymptotics,
and intensive numerical computation; however, the delicacy of these analyses/computations
has made definite conclusions elusive. One of the few definitive rules of thumb is that increasing
activation energy destabilizes detonations, while increasing overdrive or decreasing heat release
stabilizes them. However, this has been difficult to confirm globally due to difficulty/expense of
computing for sufficiently large activation energies. 
We hope that the selection of 1D results we have described indicates a clear role for the type of
dynamical systems/Evans function techniques used to study viscous shock wave, both in confirming known
rules of thumb/computational results and suggesting new possible directions of investigation--
at the same time suggesting roles for viscous theory in providing both rigor and new phenomena.

At the inviscid level, an unexpected bonus has been the discovery of the useful coordinatization of Section 
\ref{I.II}, which appears to offer useful guidance/organization of information at the level of applications.
It is to be hoped that further analysis (see open problem 3 just below) will identify similar ``master coordinates''
in the context of rNS, removing the hysteresis of Figure \ref{fig:sb}.

\medskip

{\it Open problems:}

\medskip

$\bullet$ Effects of viscosity on detonation behavior. 

$\bullet$ 1D instability of ZND detonations in the high-activation energy limit.

$\bullet$ Viscous stabilization of rNS detonations in the high-activation energy limit.

$\bullet$ Stability of weak detonations/deflagrations \cite{CF,GS} for rNS. 

\medskip

Regarding the first problem, see the interesting recent discussion by Powers and Paolucci \cite{PP} 
on complicated-chemistry reactions, pointing out that viscous length scales neglected in ZND may be on the same 
order as reaction scales important for stability.
Regarding the second, it has been addressed formally in suggestive fashion by Buckmaster--Neeves, Short, Clavin--He, etc.  
\cite{BN,S1,CH}, but up to now (a) not rigorously verified, and (b) as pointed out by
Erpenbeck, Lee--Stewart, Short, etc. \cite{Er1,LS,S1,S2}, exhibiting puzzling differences with observed numerics.
Both this and the third, hyperstabilization, problem appear to reduce to semiclassical limit/turning-point 
problems similar to those treated in Section \ref{s:II}, with governing parameter $1/\E\to 0$.
See also the related \cite{FKR} for reduced models accurately capturing behavior.
The fourth problem has been studied for simplified ``Majda''-type models in \cite{Ma1,LyZ2,Sz,LY} and for
artificial viscosity systems in \cite{LRTZ}; for a discussion in the context of the full rNS equations, see \cite{TZ4}.

\section{High-frequency stability of ZND detonations and $C^\infty$ vs. $C^\omega$ stationary phase}\label{s:II}
In this second part, we focus now on a specific topic in multidimensional stability analysis for ZND.
A delicate aspect of numerical stability investigations for ZND (inviscid) detonations is truncation of the computational domain by high-frequency asymptotics, a semiclassical limit problem for ODE.  
In this part, we focus on this issue in the most delicate {\it multi-D case}, revisiting and completing/somewhat extending
the important investigations of this topic by Erpenbeck \cite{Er3,Er4} in the 1960's.
This leads to interesting questions related to WKB expansion, turning points, and block-diagonalization/separation of modes.  In particular, as we shall describe, it highlights the difference between spectral gap and ``spectral separation,'' revealing essential differences between  $C^\infty$-coefficient and analytic-coefficient theory.  These differences are in turn related to oscillatory integrals and differences in stationary phase estimates for $C^\infty$ vs. analytic symbols.

Questions we have in mind in this section are:

\medskip

$\bullet$ Can we complete/make rigorous the turning-point investigations of Erpenbeck?

\medskip

$\bullet$ What is the meaning, finally, of such {\it inviscid} {high-frequency} results?

\subsection{Multi-d stability of ZND detonations}\label{II.I}
	The multi-D reactive Euler, or Zel'dovich--von Neumann--D\"oring (ZND) equations, in Eulerian coordinates, arr

\begin{equation}\label{ZND}
\left\{ \begin{aligned}
		\d_t \rho + \nabla_x \cdot(\rho u)  & = 0,\\
		\d_t u + \nabla_x \cdot (\rho u\otimes u)+ \nabla_x p  & = 0, \\
		\d_t E + \nabla_x \cdot (\rho u E + up)  & = 
		qk \phi(T) z ,\\
		\d_t z + \nabla_x \cdot(\rho u z) & =  - k \phi(T) z ,\\
\end{aligned}\right.
\end{equation}
where $\rho>0$ is density, $u$ velocity, 
$ E = e+ \frac{1}{2} |u|^2$ 
specific gas-dynamical energy, 
$e>0$ specific internal energy,  
and $0 \leq z \leq 1$ mass fraction of the reactant,
typically with polytropic equation of state and Arrhenius-type ignition function,
\be\label{eoseg2}
 p=\frac{\Gamma  e}{\tau}, \quad T=c^{-1} e, \quad \phi(T)=e^{\frac{-\mathcal{E}}{T}}.
\ee

\subsubsection{Planar ZND detonation waves}
\begin{figure}
\includegraphics[trim=0 50 0 50, scale=.5]{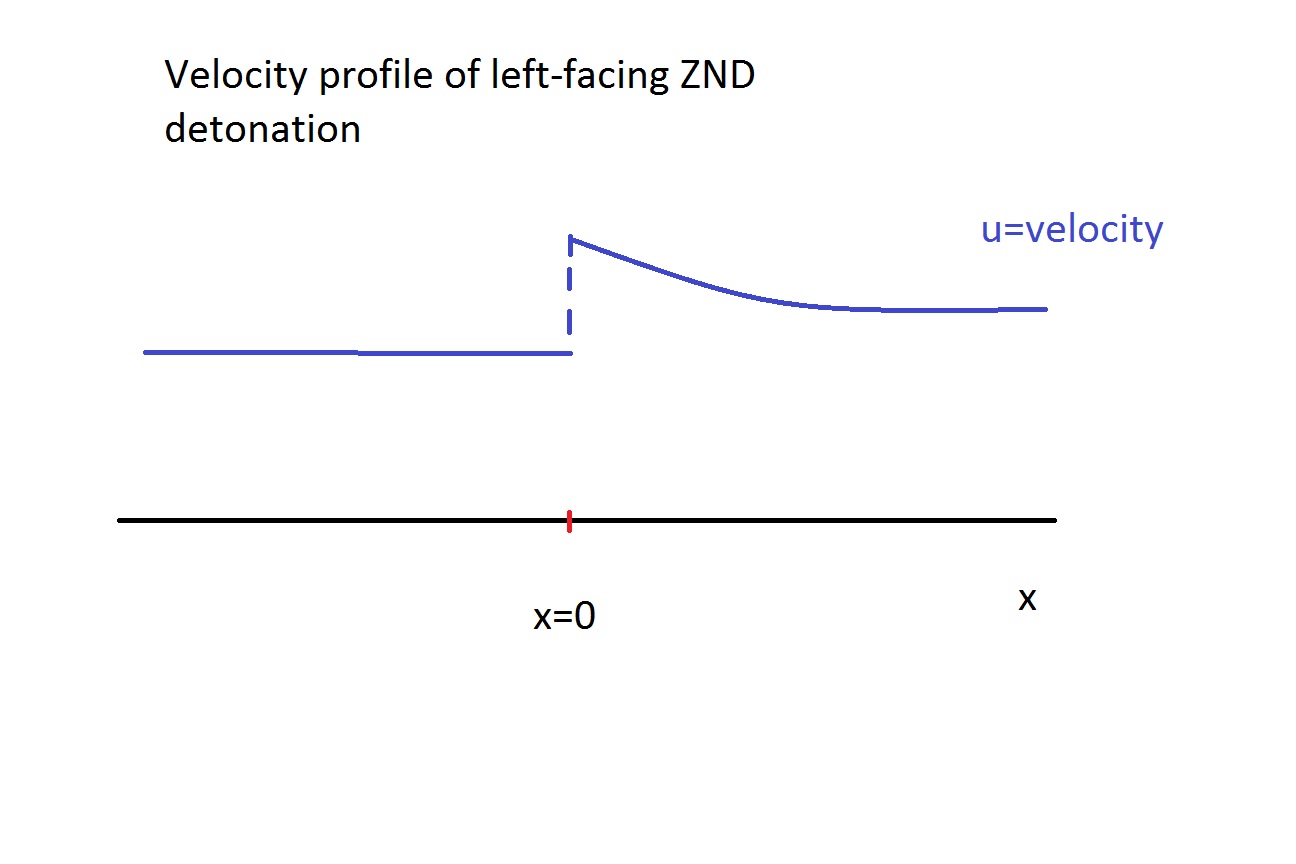}
\label{fig:multprof}
\end{figure}

A without loss of generality standing ``left-facing'' planar detonation front is a solution
\begin{equation}
(\rho, u, E, z)(x,t)=
\begin{cases}
	(\rho_-,u_-,E_-, 1),&\quad x_1<0,\\
	(\bar \rho, \bar u, \bar E, \bar z)(x_1),
	&\quad x_1\geq 0,
\end{cases}
\notag \end{equation}
of \eqref{ZND} with $(\bar \rho, \bar u, \bar E, \bar z)(x_1)\to (\rho_+,u_+,E_+, 0)$ as $x_1\to +\infty$.
This consists of a nonreactive ``Neumann'' shock at $x_1=0$, $z(0^\pm)=1$, pressurizing reactant-laden
gas moving from left to right and igniting the reaction.  
As depicted in Figure \ref{fig:multprof}, the profile is constant on $x\leq 0$ and has a reaction tail on $x_1\geq 0$, with
burned state $z= 0$ at $x_1= +\infty$.

	\subsubsection{ Spectral stability analysis}
	Consider the abstract formulation of the equations
	\be\label{abeq}
	W_t + \sum_j \partial_{x_j} F_j(W)_x=R(W). 
	\ee
Similarly as in the 1D case, a normal modes analysis leads to the 
linearized eigenvalue problem \cite{Er2,Ma2,JLW})
$ \lambda A_1^{-1} Z +  Z' \sum_{j\neq 1} i\xi_j A_j A_1^{-1} Z=EA_1^{-1} Z, $
where $Z:=A_1W$, $'$ denoting $\partial_{x_1}$, or interior equation (written without loss
of generality for simplicity in dimension $d=2$):
$$
Z'=GZ:= -(\lambda + \xi A_2 + E )A_1^{-1}Z,
\qquad x\geq 0^+,
$$
plus a (modified Rankine-Hugoniot) jump condition at $x=0$. Here and in what follows, we drop the subscript for $x$,
writing $x_1$ as simply $x$.

	\subsubsection{Evans--Lopatinski condition (Erpenbeck's stability function)}
	Normal modes $e^{\lambda t}e^{i\xi x_2}W(x_1)$, $\Re \lambda\geq 0$ correspond to zeros of the 
	Evans--Lopatinski determinant
\be\label{multilop}
D_{ZND}(\lambda)= \tilde Z_1^-(\lambda,0)\cdot
\big(\lambda[\bar W] + i\xi[F_2(\bar W)]+ R(\bar W^0)(0^+) \Big),
\ee
where $[\cdot ]$ denotes jump across $x=0$ and $\tilde Z_1$ is a (unique up to constant multiplier) 
solution of the dual equations
\be\label{multieq}
\tilde Z'=-G^* \tilde Z = (A_1^*)^{-1}(\lambda + \xi A_2 + E )^* \tilde Z 
\ee
decaying as $x\to +\infty$.  (This neat formulation, again, due to Jenssen-Lyng-Williams \cite{JLW}.)

\subsubsection{Comparison to shock wave case}\label{s:comparison}
For later, we note that the Evans--Lopatinski determinant described \eqref{multilop}--\eqref{multieq}
is quite similar to that described for the shock wave case in \cite{ZS,Z4},
with the differences that here $G$ variable-coefficient rather that constant-coefficienht and $E,R\not \equiv 0$.
In practice, \eqref{multilop} is computed numerically by approximation of \eqref{multieq} \cite{Er1,Er2,HuZ2,BZ1,BZ2}.

\subsection{High-frequency stability and the semiclassical limit}\label{II.II}
We now come to our main topic.
The numerics typically used to evaluate \eqref{multilop} are sensitive and computationally intensive,
particularly at high frequencies \cite{Er2,LS,BZ1,BZ2}; thus, an important step in obtaining reliable
results is to truncate the frequency domain by a separate, high-frequency analysis.
Even the few analytically deducible results (stability in $q\to 0$ or high overdrive limit) 
require high-frequency truncation as a crucial (and somewhat delicate) step; see Section \ref{s:I}.
Our purpose here is to:
1. Describe two recent results of Lafitte-Williams-Zumbrun on high-frequency stability \cite{LWZ1,LWZ2}: 
one instability and one stability theorem, building on the
pioneering ideas of Erpenbeck's 1960's Los Alamos Technical Report \cite{Er3,Er4}.
2. Discuss related block-diagonalization of semiclassical ODE \cite{LWZ3}.

\subsubsection{Formulation as semiclassical limit}
Setting $(\xi, \lambda)= h^{-1}(1, \zeta)$, for $|\xi|>>1$, interior equation \eqref{multieq} becomes
the semiclassical limit problem
	\be\label{semiclassical}
h\tilde Z'= (G_0 + hG_1) \tilde Z,
	\ee
where $G_0= -[(\zeta + iA_2)A_1^{-1}]]^*$ 
	involves only nonreactive gas-dynamical quantities, 
	so is identical to the symbol appearing in (nonreacting) shock 
	stability analysis, $G_1$ uniformly bounded, 
	$h=|\xi|^{-1}\to 0$.
Likewise, the boundary vector $ \big(\lambda[\bar W] + i\xi[F_2(\bar W)]+ R(\bar W^0)(0^+) \Big) $ appearing in \eqref{multilop} 
rewrites as  
\be\label{bound2}
 h^{-1}\big( \ell_0 + h\ell_1\big),
\ee
where $\ell_0=\zeta[\bar W] + i[F_2(\bar W)]$ is as in the nonreactive gas-dynamical case,
and ${R(\bar W^0)(0^+)}$ is bounded.
The difference in principle parts from the nonreactive case is just that {\it $G_0 $ is now varying in $x_1$.}

\subsubsection{Symbolic analysis}\label{symbol}
From the study of nonreactive gas dynamics \cite{Er6,Ma2,Z4,Se}, 
we know that the eigenvalues of the principal symbol $G_0$ are
\be\label{Gevals}
	\mu_1=-\kappa(\kappa\zeta+s)/\eta u_1,
\;\;\;\mu_2=-\kappa(\kappa\zeta-s)/\eta u_1,\;\;\; 
	\mu_3=\mu_4=\mu_5=\zeta/u_1,
\ee
where $\kappa=u_1/c_0$, $\eta= 1-u_1^2/c_0^2$, $c_0=$ sound speed, 
\be\label{s}
	s=\sqrt{\zeta^2+c_0^2 - u_1^2}, 
\ee
and, from the profile existence theory (specifically, the {\it Lax characteristic condition} \cite{La1,La2,Sm} on the component Neumann shock), 
\be\label{pos}
c_0^2-u_1^2>0;
\ee
here, $\mu_1$ and $\mu_2$ are {\it acoustic}, and $\mu_3,\dots,\mu_5$ {\it entropic} and {\it vorticity} modes.

Thus, on the domain $\Re \zeta  \geq 0$ relevant to the eigenvalue/stability problem, there is a single decaying mode $\mu_1$
for $\Re \zeta>0$, which extends continuously to the boundary $ \zeta=i\tau$.
For reference, we will call this the ``decaying'' mode even at points on the imaginary boundary where it becomes purely oscillatory
(as does happen for values of $\zeta=i\tau$ such that $\tau^2 \geq c_0^2-u_1^2$).
Evidently, the decaying mode $\mu_1$ remains separated from all other modes $\mu_j$ except at 
{\it glancing points} for which $\mu_1=\mu_2$, or $s=0$:
{\it equivalently, $\zeta =\pm i\sqrt{ c_0^2-u_1^2}$, a property depending on both $x_1$ and $\zeta$.}

Glancing points play a central role in the study of multi-D nonlinear stability of 
(nonreactive) viscous and inviscid shock and boundary layers \cite{Kr,Ma2,Me1,Me2,Z3,MZ,GMWZ1,GMWZ2,N},
presenting the chief technical difficulty in obtaining sharp linear resolvent bounds needed to close a nonlinear analysis.
There, the issue is to obtain bounds on a constant-coefficient symbol as frequencies $\xi, \lambda$ are varied in the
neighborhood of a glancing point.
In the present context, the problem is essentially dual: {\it for fixed frequencies $\zeta$ to understand the flow
of ODE \eqref{semiclassical} as the spatial coordinate $x$ is varied,}
a nice twist for experts in shock theory.  This leads us naturally to WKB expansion/turning point theory, 
where glancing points represent {\it nontrivial turning points}.

	\subsubsection{WKB expansion/approximate block-diagonalization}
The situation of ODE \eqref{semiclassical}, where solutions vary on a much faster scale $\sim h^{-1}$ vs. $\sim 1$
than coefficients, is precisely suited for approximation by WKB expansion.
As discussed in \cite[Section 1.1.1]{LWZ3}, WKB expansion is closely related to the {\it method of 
repeated diagonalization} \cite{Lev,CL}, both methods consisting of constructing approximate solutions from 
diagonal modes of a sufficiently high-order decoupled system. 

\medskip
{\it Primitive version:} We illustrate the approach by a treatment of the simplest (nonglancing) case, when $\mu_1$ and
$\mu_2$ remain separated for all $x\geq 0$.
This occurs, for example, on the strictly unstable set $\Re \zeta>0$.
Then, the decaying mode $\mu_1$ remains separated from the remaining eigenvalues $\mu_1,\dots, \mu_5$ of $G_0(x)$.
By standard matrix perturbation theory \cite{K}, it follows that there exists a change of coordinates $T_{G_0}$, depending
smoothly on $G_0$, such that
\be\label{diag0}
T_{G_0}^{-1}G_0T_{G_0} =\bp \mu_1 & 0 \\ 0 & M_{G_0}\ep.
\ee

Setting $T(x):=T_{G_0(x)}$, $M(x):=M_{G_0(x)}$, and making the change of coordinates $\tilde Z(x)=T(x)\tilde W(x)$,
we convert \eqref{semiclassical} to an ODE 
\be\label{intode}
h\tilde W'= \bp \mu_1 & 0 \\ 0 & M\ep \tilde W-  (hT^{-1}T'\tilde ) W,
\ee
that, to order $O(h)$ of the commutator term $hT^{-1}T'$ , is {\it block-diagonal with a decoupled $\mu_1$ block.}

Next, observe that an $O(h)$ perturbation of a block-diagonal matrix with spectrally separated blocks may be block-diagonalized by a coordinate
change $T_2=\Id +O(h)$ that is a smooth $O(h)$ perturbation of the identity \cite{K}; applying such a coordinate change, and observing that the associated
commutator term $hT_2^{-1}T_2'= hT_2^{-1}(O(h))'= O(h^2)$, we can thus reduce to an equation that is block-diagonal to $O(h^2)$.
Repeating this process, we may obtain an equation that is {\it block-diagonal up to arbitrarily high order error $O(h^p)$}, so long as 
the coefficients of the original equation \eqref{semiclassical} possess sufficient regularity that derivatives in commutator terms remain $O(1)$. 

Untangling coordinate changes, this  suggests that the unique solution $\tilde Z_1$ decaying as $x\to +\infty$
	``tracks'' to $O(h)$ the $R_1$ eigendirection associated with $\mu_1$, satisfying the WKB-like approximation
$$
\tilde Z_1(x)= e^{h^{-1}\int_0^{x} (\mu_1+O(h))(y)dy} (R_1(x)+O(h)),
$$
with in particular $\tilde Z_1(0)=R_1+O(h)$, where $R_1$ is an eigenvector of the decaying mode of $-G_0^*(0^+)$.

Recall \cite{Er6,Ma2,ZS} that the Lopatinski determinant for the component Neumann shock is 
\be\label{lopform}
D_N= \ell_0\cdot R_1,
\ee
where $\ell_0$ is the principal part of boundary vector \eqref{bound2}.
Thus, {\it assuming that the above approximate diagonalization procedure with formal error $O(h^p)$ may
be converted to an exact block-diagonalization with rigorous convergence error $O(h^p)$,}
we may conclude thag
\be\label{ZND_N}
D_{ZND}(\xi, \lambda)=D_{N}(\xi,\lambda) (1+O(h)),
\ee
where $D_{N}$ is the Lopatinski determinant for the stability problem associated with the Neumann shock at $x=0$, 
{\it hence ZND detonation is high-frequency stable for such choices of $\zeta$}
(which include always the strictly unstable set $\Re \zeta>0$) 
{\it if and only if its component Neumann shock is stable.}

\medskip

{\it The glancing case.} In the glancing case, $s(x_*,\zeta_*)=0$ for some $x_*\geq 0$, and there is a nontrivial turning point
at $x=x_*$.
In this case, for $\zeta$ and $x$ local to $\zeta_*, x_*$, there is no uniform separation between $\mu_1$ and $\mu_2$,
and the above-described complete diagonalization procedure no longer works.
However, observing that $\mu_1$ and $\mu_2$ together remain spectrally separated from $\mu_3,\dots, \mu_5$,
we can still approximately block-diagonalize to a system with coefficient $\bp P& 0\\ 0 & N\ep$, where
$P$ is a $2\times 2$ block corresponding to the total eigenspace of $G_0$ associated with $\mu_1$ and $\mu_2$,
in particular having eigenvalues $\mu_1$ and $\mu_2$.
It is shown by a normal form analysis in \cite{LWZ2} that any such $2\times 2$ block, under a nondegeneracy condition
on the variation of its eigenvalues with respect to $x$ at $x_*$, can be reduced further to an arbitrarily
high-order perturbation of {\it Airy's equation}, written as a $2\times 2$ system, where in this case the nondegeneracy
condition is just
\be\label{nondeg2}
(s^2)'(x_*)\neq 0.
\ee
{\it Assuming as before that the above approximate diagonalization procedure 
be converted to an exact block-diagonalization with rigorous convergence error,}
we may thus hope to analyze this case by reference to the known (see, e.g.: \cite{AS}) behavior of the Airy equation.

	\subsection{The Erpenbeck high-frequency stability theorems}\label{II.III}
We are now ready to state our main theorems regarding profiles of the abstract system \eqref{abeq}.
We make the following assumptions:

\begin{assumption}\label{hypass}
The associated nonreactive system $W_t + \sum_j \partial_{x_j} F_j(W)_x=0$ is hyperbolic for all
value of $W$ lying on the detonation profile $\bar W(x)$.
\end{assumption}

\begin{assumption}\label{lopass}
The component Neumann shock for:profile $\bar W$ is Lopatinski stable.
	\end{assumption}

	\begin{assumption} \label{analyticass}
		The coefficients of system \eqref{abeq} are {\it real analytic}.
	\end{assumption}

	\begin{definition} \label{typeass}
		A detonation is type $I$ (resp. $D$) if $c_0^2-u_1^2$ is increasing (resp. decreasing).
	\end{definition}

	\br\label{classrmk}
	Erpenbeck classifies a number of materials/detonations as class I or D.  
	More general cases may in principle be treated by elaboration of the techniques
	used to treat classes I and D.
	\er

	\bt[LWZ2012]\label{t1}
	Under Assumptions \ref{hypass}--\ref{analyticass}, plus an additional (frequently satisfied) ratio condition, 
	type I detonations exhibit Lopatinski instabilities of arbitrarily high frequency.
	\et

\begin{proof}
[Sketch of Proof] (case of turning point) 
By the block-diagonalization procedure described above, first reduce to a $2\times 2$ block involving only the growth
modes $\mu_1$ and $\mu_2$.
For type I, growth rates $\mu_1$ and $\mu_2$ correspond
to exponentially growing/decaying modes for $x_1>x_*$, 
oscillatory modes for $x_1<x_*$, the connections between these solutions across the value
$x=x_*$ being determined by behavior of the Airy equation.
The question is whether the Airy equation takes the pure decay mode to the corresponding pure oscillatory mode 
(the ``decaying'' mode at $x=0^+$).  It does not-- rather to the average of the two decaying modes \cite{AS}, 
giving a solution composed of oscillating comparable-size parts, which, under the ratio condition, 
cancel for a lattice of $ \lambda=h^{-1}i\tau +\nu$ with $\Re \nu>0$.  
(Otherwise they cancel for frequencies $\Re \nu<0$ not giving instability.) 
\end{proof}

\bt[LWZ2015]\label{t2}
Under Assumptions A1-A3, type D detonations are Lopatinski stable for sufficiently high frequencies.
\et

\begin{proof}[Sketch of Proof] (case of turning point) 
	As in case I, the problem reduces to a $2\times 2$ block, and the study of connections 
	across the turning point $x=x_*$ determined by behavior of the Airy equation.
	For type D, the reverse happens, 
	By reflection symmetry of the Airy equation, there holds in case D essentially the reverse situation to
	that of case I, featuring oscillatory modes for $x>x_*$ and exponentially growing/decaying modes for
	$x<x_*$, connected by a reverse Airy flow.
	So, again we see that the pure
	``decay'' (now actually oscillating) mode at $+\infty$ does not connect to the pure growth mode at $x=x_*^-$,
	but contains at least some component of the (actual) ``decay'' mode for $x<x_*$.
	It follows by order $e^{O(x/h)}$ exponential growth in the backward direction of this decaying mode,
	together with order $e^{O(-x/h)}$ exponential decay in the backward direction of the complementary
	growing mode, that the solution at $x=0$ is dominated by the decay-mode component
	Thus, $\tilde Z(0^+)$ lies to exponentially small order in the $R_1$ direction, $R_1$ 
	as in \eqref{lopform},
	giving the (stable) shock Lopatinski determinant in the limit, as in the simplest (nonglancing) case.
\end{proof}

{\it Technical issues:}
	1. Exact vs. approximate block-diagonalization.
	2. Block-diagonalization at $+\infty$.
	3. Turning points at $x_*=0, +\infty$, and exact vs. approximate conjugation to Airy/Bessel
	(daunting).

	\begin{remarks}\label{erprmks}
		1. Theorems \ref{t1}-\ref{t2} justify the voluminous literature on numerical 
		multi-d stability stability results, full implications of which were previously unclear.

		2. The arguments streamline/modernize the analysis of \cite{Er3,Er4} 
		(carried out originally by WKB expansion in all $5$ modes!).
	But also new analysis at degenerate frequencies is needed for the complete stability result.

	3. The proofs are still hard work! (Amazing achievement of Erpenbeck in the 1960's.

	4. We have suppressed discussion of conjugations to Airy/Bessel equations (difficult! the latter new),
	and the related huge contributions of F. Olver in asymptotics of special functions \cite{O,AS}.
\end{remarks}

\subsection{Exact block-diagonalization and $C^\infty$ vs. $C^\omega$ stationary phase}\label{II.IV}
	Consider an approximately block-diagonal equation 
	$$hW'=AW+h^p\Theta,$$
	$\Theta=$ error, and seek $T=\bp I& h^p\alpha_{12}\\h^p \alpha_{21}& I\ep$ 
	such that $W=TZ$ gives $hZ'=DZ$, exact.
Equating first diagonal, then off-diagonal blocks in
$ (hT'+TD)Z= (A+h^p\Theta)TZ, $ yields Ricatti equations
\be\label{ric}
\begin{aligned}
	h \alpha_{12}' &= 
		A_{11}\alpha_{12}- \alpha_{12}A_{22}  
	+  
	\Theta_{12}- h^{2p}\alpha_{12} \Theta_{21}\alpha_{21}- h^p \Theta_{11} \alpha_{21} ,\\
	h \alpha_{21}' &= 
	A_{22}\alpha_{21}- \alpha_{21}A_{11}  +  
	\Theta_{21}- h^{2p}\alpha_{21} \Theta_{12}\alpha_{12}- h^p \Theta_{22} \alpha_{21}  ,\\
	\end{aligned}
\ee
or, viewed as a block vector equation in $\alpha=(\alpha_{12},\alpha_{21})$:
\be\label{vecversion}
h\alpha'= \cA(0) \alpha 
+
(\mathcal{A}(z)-\mathcal{A}(0))\alpha +  Q(\alpha,\Theta,h)
.
\ee

{\it Observation} Sylvester equation, hence $\sigma(A_{11}) \cap \sigma(A_{22})=\emptyset$ implies $0\not \in \sigma(\cA(0))$.  

	\subsubsection{Lyapunov-Perron formulation (standard)}
\begin{figure}
\includegraphics[trim=0 40 0 40, scale=.5]{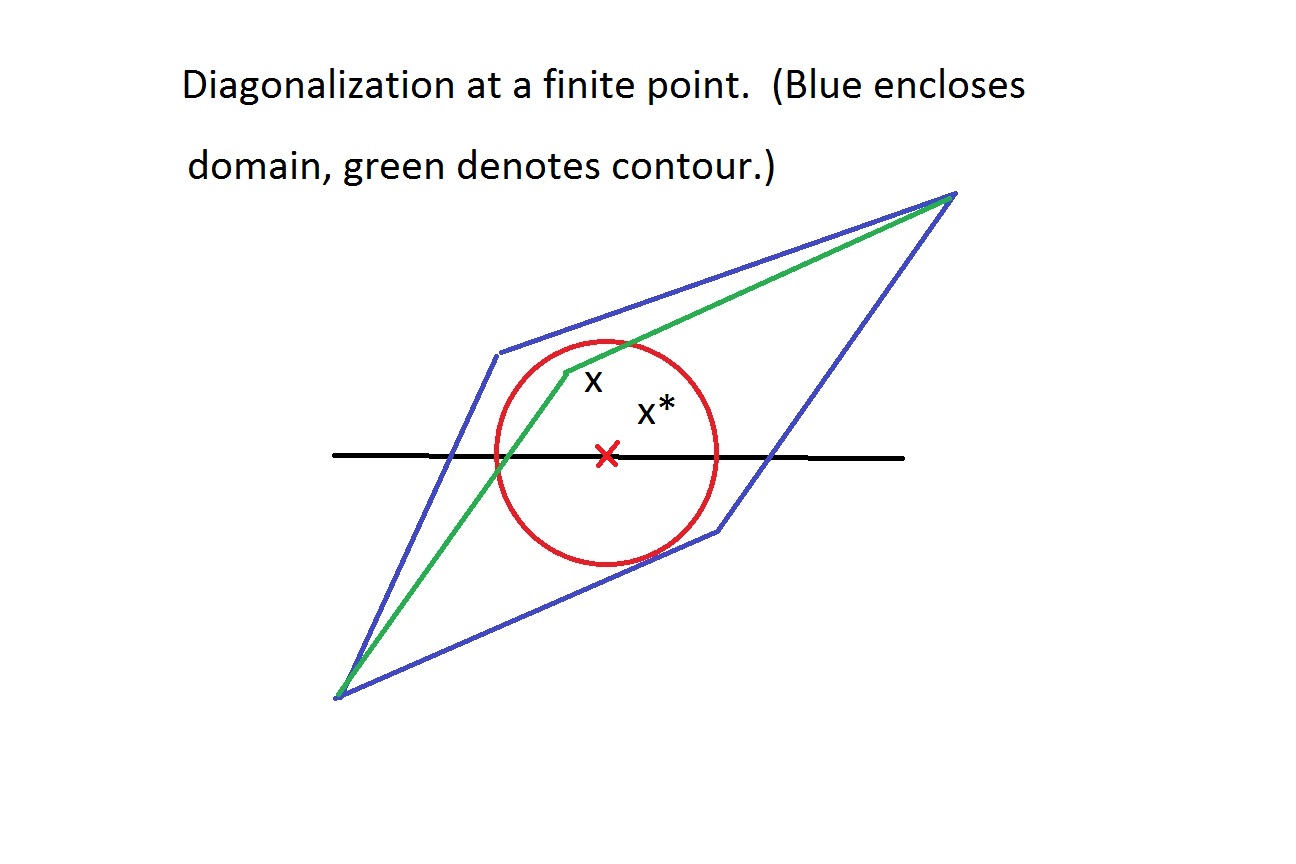}
\caption{Block diagonalization at a finite point.}
\label{fig:finite}
\end{figure}

	From $  h\alpha'= \cA(0) \alpha  +
(\mathcal{A}(z)-\mathcal{A}(0))\alpha +  Q(\alpha,\Theta,h)$,
we obtain by Duhamel's principle the integral fixed-point equation 
$$
\begin{aligned}
\alpha(x)=  
\CalT \alpha(x)&:=
  h^{-1}\int_{z_*}^x e^{h^{-1}\cA(0) (x-y)}
  \Pi_U \big( (\mathcal{A}(y)- \mathcal{A}(0))\alpha(y) + Q
  (y) \big) \,dy \\
 &\quad +  h^{-1}\int_{z^*}^x e^{h^{-1}\cA(0) (x-y)}
  \Pi_S \big( (\mathcal{A}(y)- \mathcal{A}(0))\alpha(y) + Q
  (y) \big) \,dy 
,
\end{aligned}
	$$
on diamond 
$ \mathcal{D}:=\{ x:\, |\arg \big((x-z_*)/\gamma\big)|, \; |\arg \big( (z^*-x)/\gamma \big)|\leq \eps \} $,
where $\gamma \in \C$, $|\gamma|=1$ is chosen so that 
$\cA(0) \gamma$ has spectral gap, 
and $\Pi_U$, $\Pi_S$ denote stable/unstable projectors of $\cA(0)\gamma$; see Figure \ref{fig:finite}.
 Mapping $\CalT$ is contractive by $O(e^{-\eta |x-y|/h})$ decay of propagators, plus smallness of the source.

 \br\label{escapermk}
Here, we have used analyticity to escape the real axis and recover a spectral gap.
This is essentially a finite-regularity version of a theorem of Wasow \cite{W}
in the $h$-analytic case \cite{LWZ3}.
\er

\subsubsection{Block diagonalization at infinity}
In many problems (e.g., detonation), we must treat unbounded intervals, diagonalization at infinity.
This may be carried out by the following modifications of the argument for
the finite-turning point case \cite{LWZ3}. Briefly, we:

\smallskip

$\bullet$ Require analyticity on a wedge about infinity, not just neighborhood of real axis
(and verify that this is indeed guaranteed by stable manifold construction for analytic coefficient profile equation).

\smallskip

\begin{figure}
\includegraphics[trim=0 40 0 40, scale=.5]{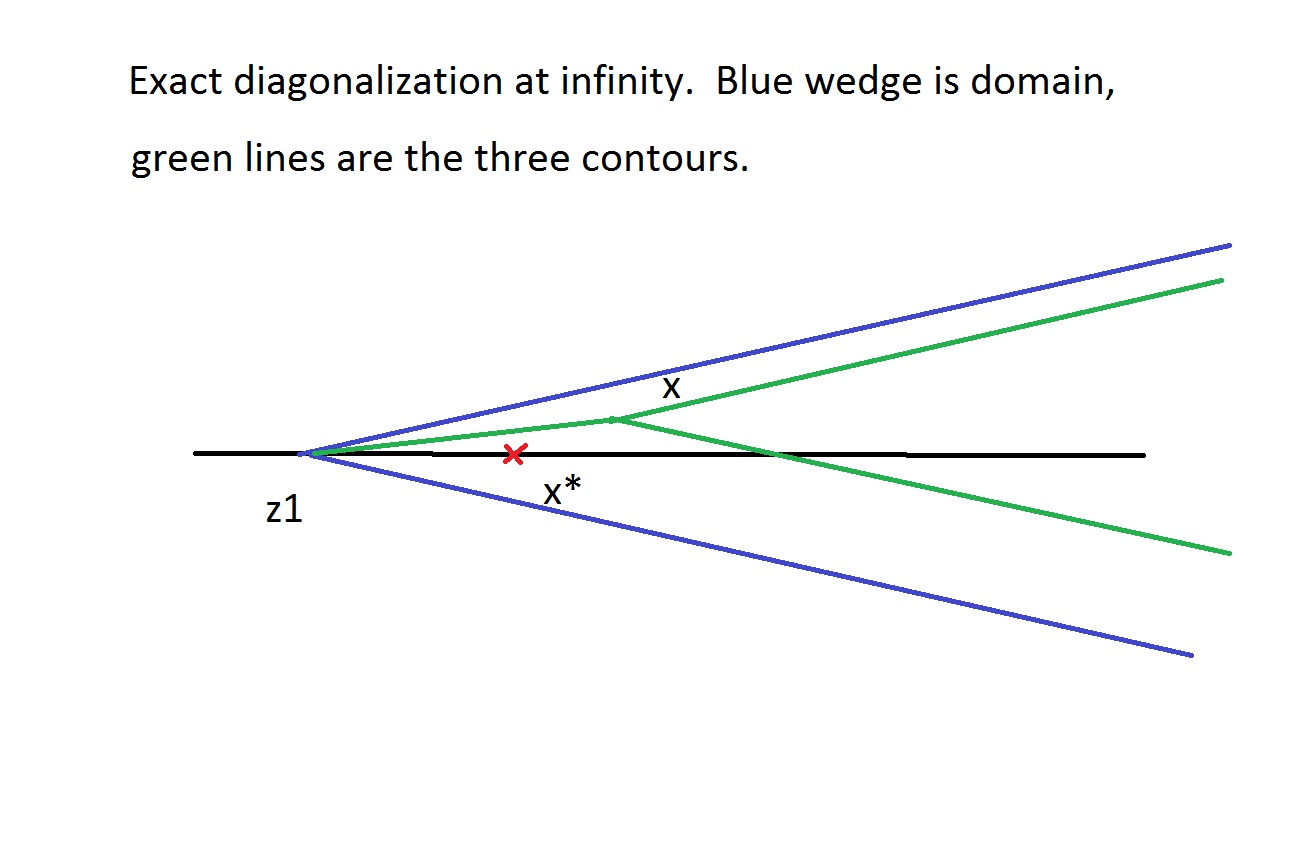}
\caption{Block diagonalization at infinity: contour configuration.}
\label{fig:infty}
\end{figure}

$\bullet$ Use three contour directions to recover a spectral gap; see Figure \ref{fig:infty}.

\bigskip

	\subsection{Counterexamples and $C^\infty$ vs. $C^\omega$ stationary phase}\label{II.V}
	Our treatment of multi-D high-frequency behavior has a different flavor from the 
	analyses of 1-D stability in part one: in particular, we have used analyticity of coefficients 
	and moved away from the real line.

	\medskip

	{\it A natural question:} Is this necessary?  In particular, could we by some other method perform 
	block-diagonalization for $C^\infty$ (or just $C^r$ as in the 1-D case) coefficients?  

\medskip
	
{\it Rephrasing:}
1. Is a spectral gap between blocks (as in classical ODE techniques \cite{L}) 
needed for exact $C^\infty$ diagonalization, or just spectral separation?
	2. And (Wasow, 1980's \cite{W}), can analytic block-diagonalization be performed 
	globally under appropriate global assumptions?

	\subsubsection{Counterexamples: reduction to oscillatory integral}
	The answer to both of the above questions is ``no'' \cite{LWZ3}, as we now describe.
Consider the $2\times 2$ triangular system
\ba\label{ntrisys}
hW'=
\mathcal{A}(x,h)W:=
\bp \lambda_1(x)& h^p \theta(x)\\0 & \lambda_2(x)\ep W , \qquad W\in \C^2,\; p\geq 1,
\ea
$\theta$ uniformly bounded, with globally separated eigenvalues $\lambda_1(x)= x+i$, $\lambda_2=-(x+i)$.


\bl[\cite{LWZ3}]\label{l:oscequiv}
There exists $T(x,h)$ on $[-L,L]\subset \R$, $0\leq h\leq h_0$,
$T,$ $T^{-1}$ uniformly bounded in $C^1$, for which 
$W=TZ$ converts \eqref{ntrisys} to a diagonal system $hZ'=D(x,h) Z$, if and only if
\be\label{e:onlyif}
\hbox{\rm $\int_{-x}^x e^{-y^2/h- 2iy/h} \theta(y) dy \lesssim h e^{-x^2/h}$
for all $|x|\leq L$.  }
\ee
\el

\begin{proof}[Sketch of proof]
	It is sufficient to seek a triangular diagonalizer $T=\bp 1 & h^p\alpha\\ 0 & 1\ep$, in which case
	Ricatti equation \eqref{ric} reduces to a {\it scalar, linear ODE} in $\alpha$:
	\be\label{scalareq}
		h \alpha'= (\lambda_1-\lambda_2)\alpha +  \theta .
		\ee
By Duhamel's principle/variation of constants, existence of a uniformly bounded $T$ thus implies
uniform boundedness of
$$
\begin{aligned}
\alpha(x,h) &=
h^{-1}\int_0^x e^{h^{-1}\int_y^x (\lambda_1 -\lambda_2)(z)dz} \theta(y) dy
- e^{h^{-1}\int_0^x (\lambda_1 -\lambda_2)(z)dz} \alpha(0,h) \\
&=
e^{(x^2+2ix)/h}\Big(h^{-1}\int_0^x e^{ -(y^2+2iy)/h} \theta(y) dy - \alpha(0,h)\Big).
\end{aligned}
$$

The direction {\it ($\Leftarrow$)} then follows by $ e^{-2ix}\alpha(x) - e^{2ix}\alpha(-x) =
e^{x^2/h}\int_{-x}^x e^{ -(y^2+2iy)/h} \theta(y) dy. $

The direction {\it ($\Rightarrow$)} follows by direct computation, choosing
$\alpha(0,h)= h^{-1}\int_0^L e^{ -(y^2+2iy)/h} \theta(y) dy $. 
\end{proof}

	\subsubsection{Failure of global conjugators}

	\bl[\cite{LWZ3}]\label{quadstat}
For $a\not \equiv 0$ analytic on $[-L,L]\times [-i,i]$, and $h\to 0^+$, 
\be\label{e:quadtrans}
\int_{-x}^x e^{-y^2/h- 2iy/h} a(y) dy 
\begin{cases}
	\lesssim h e^{-\frac{x^2}{h}}, & 0<x\leq L < 1,\\
	\sim h^{(j+1)/2} e^{-\frac{1}{h}}, & 1< c_0\leq x\leq L,\\
\end{cases}
\ee
where $j=$ order of first nonvanishing derivative of $a$ at $z=i$.
\el

\begin{proof}
	The general case follows by complex-analytic stationary phase estimates (see \cite{M,PW}). 
	The simplest case $a\equiv 1$ (enough for a counterexample), follows from 
	$\int_{-\infty}^{+\infty} e^{-y^2/h- 2iy/h} a(y) dy =e^{-1/h}$, which
	follows from the fact that the Fourier transform of a Gaussian is 
	Gaussian.
\end{proof}

{\it Consequence:}
	Lemma \ref{quadstat} implies that there is no bounded block-diagonalizing conjugator of
	\eqref{ntrisys} on $[-x,x]$ for $|x|>1$, resolving a 30-year open question of Wasow \cite{W}.

	\subsubsection{Failure of local conjugators for $C^\infty$ coefficients}

\bl[\cite{LWZ3}]\label{halfprop}
For $0<c_0\leq x$ and
$a(y):= e^{-y^{-1/(s-1)}}$ for $y>0$ and $0$ for $y\leq 0$,
$$
\int_{-x}^x e^{-y^2/h-2iy/h} a(y) dy \sim
h^{1-1/2s} e^{\frac{-c(s)+ d(s) h^{1-1/s} + O(h^{2(1-1/s)}) }  {h^{1/s}} } ,
$$
$1<s<\infty$, as $h\to 0^+$, where $c(s)>0$, and $\Re d(s)
<0$ for $s<2$.
\el

\begin{proof}[Sketch of proof]
Defining $\alpha= 1-1/s \in (0,1)$, $\beta= e^{-\frac{i\pi (1-1/s)}{2} }$,
	deform contour $[0,+\infty]$ 
	to
	$ z=h^\alpha \beta t$, $t\in (0,+\infty)$,
	to obtain
	$
	I(h)\sim h^\alpha \int_0^\infty e^{\frac{i\beta( -2t -t^{-\theta} +
	i\beta h^\alpha  t^2)} {h^{1-\alpha}}  }dt , 
	$
	then apply a standard stationary phase estimate about the
	nondegenerate maximum of the phase at $t_0=  2^{-(1-1/s)} + O(h)$.
\end{proof}

\medskip
	{\it Moral:} Results may vary for $C^\infty$ coefficients!

\medskip
{\it Related phenomena:}
	1. Diffraction by $C^\infty$ vs. analytic boundary in $\R^3$
	\cite{Leb}.  2. Probability one of Weyl distribution (``cloud'') for asymptotic spectrum of a random 
	$C^\infty$-perturbation of an analytic-coefficient operator with asymptotic:spectra initially confined
	to a curve \cite{HS}.

\subsubsection{Coda: Gevrey-regularity stationary phase}\label{gevstat}
For Gevrey norm $ \| a\|_{s,T}:= \sup_{j}|\partial_x^j a| (j!)^s/T^{j} $,
define the Gevrey class $\mathcal{G}^{s,T}$ of functions with bounded Gevrey norm.
Here, $s=1$ corresponds to analyticity on a strip of width $T$ about the real axis $\R$, 
while $s\to \infty$ corresponds to absence of regularity, with Gevrey-class functions interpolating between.
The following result gives an upper bound corresponding to the lower bound of Lemma \ref{halfprop}.

\bpr[\cite{LWZ3}]\label{p:gevrey}
For $a\in \mathcal{G}^{s,T_0}$ on $[-L,L]$, $T_0, T>1$, 
$|x|\leq L$, and some $c=c(T_1,T,s)>0$,
\be\label{comp}
\int_{-x}^x e^{-y^2/h- 2iy/h} a(y) dy \lesssim h^{1/2} \|a\|_{T,s}  e^{-c/h^{1/s}}.
\ee
\epr

Proposition \ref{p:gevrey} interpolates between the algebraic $O(h^{r})$ van der Korput bounds for $C^r$ symbols 
(roughly, $s=\infty$) and the exponential $O(h^{1/2}e^{-1/h})$ bounds for analytic symbols $a$ obtained by the saddlepoint method/analytic stationary phase.  Lemma \ref{halfprop} shows that \eqref{comp} is sharp.

\medskip

(Proof by Fourier cutoff/standard complex-analytic stationary phase.)

\subsection{Discussion and open problems}\label{disc2}
Our turning-point analyses in the first part of this section completes and somewhat simplifies
the high-frequency stability program laid out by Erpenbeck in the 1960's,
in his tour de force analyses \cite{Er3,Er4}.
This in turn solidifies the foundation of the many (and delicate)
numerical multi-D stability studies for ZND, by rigorously truncating the computational frequency domain.
On the other hand, our analysis in the second part of this section on sensitivity of block-diagonalization/WKB expansion
with respect to $C^\infty $ (indeed, Gevrey-class) perturbations raises interesting philosophical questions
about the physical meaning of our multi-D high-frequency stability results, as intuitively we think of physical coefficients
as inexactly known.

Recall that the 1-D high-frequency stability results of \cite{Z1} used a different, $C^r$ diagonalization method,
so this issue does not arise in 1-D.
Likewise, smooth dependence on coefficients with respect to $C^r$ perturbation of the Evans-Lopatinski determinant 
$D_{ZND}(\xi,\lambda)$ restricted to compact frequency domains \cite{PZ} implies 
that the {\it strictly unstabilities} asserted for analytic coefficients in Theorem \ref{t1} persist under $C^r$ 
perturbations of the coefficients, so there is no issue for our instability results.
That is, the Evans function is itself robust, independent of the methods that wr used to estimate it.
Even in the stable case, we obtain from this point of view robust stability estimates on any bounded domain, no matter how
large, in particular for domains far out of practical computation range.  Thus, the results of Theorem \ref{t2} 
have practical relevance in this restricted sense independent of questions regarding analyticity of coefficients.
The philosophical resolution of the remaining issue for ultra-high frequencies, may perhaps, similarly as other issues touched on
in Section \ref{s:I}, lie in the inclusion of transport (viscosity/heat conduction/diffusion) effects, which stabilize
spectrum for frequencies on the order of one over the size of associated coefficients.


\medskip

{\it Open problems:}

\medskip

	$\bullet$ ZND limit for multi-d (interaction of viscosity, turning points).

	\medskip

	$\bullet$ Multi-d numerics for rNS (no apparent obstacle, but computationally intensive).

	\medskip

	$\bullet$ Rigorous analysis of 1-d viscous hyperstabilization (again, apparent interaction
	of turning points vs. viscous effects).

	\medskip



\end{document}